\documentclass[a4paper,11pt]{article}
\usepackage{graphicx}
\usepackage{amsfonts}
\usepackage{amsmath}
\usepackage[utf8]{inputenc}
\usepackage{url}
\usepackage{color}

\newtheorem {Proposition}{Proposition}[section]
\newtheorem {Lemma}[Proposition] {Lemma}
\newtheorem {Theorem}[Proposition]{Theorem}
\newtheorem {Corollary}[Proposition]{Corollary}
\newtheorem {Remark}[Proposition]{Remark}

\numberwithin{equation}{section}

\title{Berry-Esseen bounds for weighted averages of Poisson avoidance functionals}
\author{Eustasio del Barrio\footnote{Research partially supported by
grant MTM2011-28657-C02-01 of the Spanish Ministerio de Educacion y Ciencia and by the Consejer\'{\i}a de Educaci\'on de la 
Junta de Castilla y Le\'on,
grant VA212U13}\\
{\it Universidad de Valladolid, IMUVA}}
\date{May 4, 2015}

\begin{document}

\maketitle
\begin{abstract}
We consider functionals which are weighted averages of the avoidance function of a Poisson process. 
Using the approach to Stein's method based on Malliavin calculus for Poisson functionals 
we provide explicit bounds for the Wasserstein distance between these standardized 
functionals and the standard normal distribution. Our approach relies on closed-form expressions for the action of 
some Malliavin type operators on avoidance functionals of Poisson  processes. As a 
result we provide Berry-Esseen bounds in the CLT for the volume of the union of balls of a fixed radius 
around random Poisson centers or for the quantization error around points of a Poisson process.
We also give Berry-Esseen bounds for avoidance functionals of empirical measures.
\end{abstract}

\section{Introduction}

Stein's method of normal approximation has become one of the main tools
for proving central lmit theorems for general functions of independent
random variables. In the last years the power of the method has been 
greatly expanded in a number of directions. In a series of papers 
(\cite{NourdinPeccati2009a}, \cite{NourdinPeccati2009b}, \cite{NourdinPeccatiReveillac2010})
Nourdin, Peccati and co-authors showed how the ideas of Stein's method could be combined
with Malliavin calculus to produce explicit bounds for the normal approximation
of smooth Gaussian functionals. More recently, in \cite{Peccatietal2010}, the method was 
extended to cover normal approximation in Wasserstein distance of functionals of 
Poisson random measures. This approach has been successfully used for proving
CLT's for sequences of multiple Wiener-It\^o integrals with respect to a Poisson measure
(see \cite{Peccatietal2010}, \cite{ReitznerSchulte2013} and \cite{Lastetal2014}) and in several instances
in stochastic geometry, including CLT's for the volume of the Poisson-Voronoi 
approximation of a compact convex set in Euclidean space (see \cite{Schulte2012}),
for the length of a random geometric graph (see \cite{ReitznerSchulte2013})
or for geometric functionals of intersection processes of a Poisson process of $k$-dimensional flats in 
$\mathbb{R}^d$ (see \cite{Lastetal2014}).

The key to the Malliavin calculus approach to Stein's method for Poisson functionals is that the
Wasserstein distance between a standardized Poisson functional and the standard normal
distribution is bounded by moments of some Malliavin operators acting on the functional. We 
refer to \cite{LastPenrose2011} and \cite{Peccatietal2010} for a complete account
of the theory and to Sections 2 and 3 below for a succint, self-contained description of the
main facts about it. Now, if  $\eta_\lambda$ is a homogeneous Poisson point process on $\mathbb{R}^d$ with intensity $\lambda>0$, 
$F_\lambda=F_\lambda(\eta_\lambda)$ satisfies $0<\mbox{Var}(F_\lambda)<\infty$ and
$X$  is a standard normal random variable then 
\begin{equation}\label{cota1}
d_W\left(\frac{F_\lambda-\mbox{E}(F_\lambda)}{\sqrt{\mbox{Var}(F_\lambda)}},X \right)\leq 
\frac{(\mbox{Var} (\langle D F_\lambda, -D L^{-1}(F_\lambda-\mbox{E}(F_\lambda))\rangle))^{1/2}+({\textstyle \int 
\mbox{E}(D_z F_\lambda)^4 d\mu(z) })^{1/2}}{\mbox{Var}(F_\lambda)},
\end{equation}
where $D$ denotes the diference operator and $L$ the Orstein-Uhlenbeck operator. 
More recently, in \cite{EiselbacherThale2013}, a similar, more involved  bound was proved for the case
of $d_K$, the Kolmogorov metric.
Both 
$D F_\lambda$ and $-D L^{-1}(F_\lambda-\mbox{E}(F_\lambda))$ have a chaos decomposition that can be simply
expressed in terms of the chaos decomposition of $F_\lambda$, namely, in the expansion
\begin{equation}\label{PoissonChaos0}
F_\lambda(\eta_\lambda)=\mbox{E}(F_\lambda(\eta_\lambda)) +\sum_{n=1}^\infty I_n(f_{n,\lambda}),
\end{equation}
where $I_n(f_{n,\lambda})$ denotes the multiple Wiener-It\^o integral of $f_{n,\lambda}$ with respect to 
$\eta_\lambda$ (we refer again to refer to \cite{LastPenrose2011}, \cite{Peccatietal2010} and Section 2 for 
details). This enables one to control the upper bound in (\ref{cota1}) through the control of moments
of products of multiple Wiener integrals of different orders. When the chaos expansion (\ref{PoissonChaos0})
consists of a finite number of terms this can be done through
the use of diagram formulae, as discussed in \cite{Peccatietal2010}, \cite{ReitznerSchulte2013}
or \cite{Lastetal2014}. However, the kind of bound that one obtains from these diagram formulae is 
not tight enough for functionals with infinite terms in their chaos expansion. Still, the approach
can be adapted through truncation arguments, as in \cite{Schulte2012} or \cite{Lastetal2014} to prove
CLT's, but, to our best knowledge, no Berry-Esseen bounds can be derived from this type of approach
in the case of infinite chaos expasions.

A more direct approach to obtain upper bounds for the right-hand side in (\ref{cota1}) would
be to exploit \textit{pathwise} representations of the operators $D$ and $D L^{-1}$, not relying
on $L_2$ expansions. There is indeed a simple pathwise representation for $D$ 
(see (\ref{differenceoperator}) below). A major breakthrough in the theory is given by \cite{LastPeccatiSchulte2014},
which provides a simple representation of the action of $L^{-1}$ through a so-called \textit{Mehler's formula}
(see Theorem 3.2 in \cite{LastPeccatiSchulte2014}). This representation greatly
increases the usability of upper bounds like (\ref{cota1}). In the cited reference, 
the represantation is combined with a Poincaré inequality (see Proposition
2.5 in \cite{LastPeccatiSchulte2014} to produce simple bounds for the distance ($d_W$ or $d_K$)
to normality of Poisson functionals (Theorems 1.1 and 1.2 in the cited reference).

In this paper we consider a special class of functionals for which the target
functionals admit a 
particularly simple
pathwise representation. We will assume throughout the paper that $\eta_\lambda$ is a
homogeneous Poisson point process on $\mathbb{R}^d$ with intensity $\lambda>0$, defined on a probability space 
$(\Omega,\mathcal{F},P)$, $\mathcal{A}$ a bounded open set on $\mathbb{R}^d$ and $(\mathcal{B},\mathcal{G},\nu)$
a measure space. We further assume that, for each $s\in\mathcal{B}$, $Q(s)$ is a Borel set on $\mathbb{R}^d$
chosen in such a way that the map $(\omega,x,s)\mapsto \mathbf{1}(\eta_\lambda (x+\lambda^{-1/d}Q(s))=0)$
is $\mathcal{F}\otimes \beta^d \otimes \mathcal{G}$ measurable, where $\beta^d$
denotes the Borel $\sigma$-field on $\mathbb{R}^d$ and $\mathbf{1}$ stands for the indicator of a set.
Under these assumptions, we will focus on the functional
\begin{equation}\label{functionaldefinition}
F_\lambda(\eta_\lambda):=\int_{\mathcal{A}\times\mathcal{B}} \mathbf{1}(\eta_\lambda (x+\lambda^{-1/d}Q(s))=0) d(\ell\otimes \nu)(x,s).
\end{equation}
This functional is a weighted average of values of the \textit{avoidance functional} $A\mapsto \mathbf{1}(\eta_\lambda(A)=0)$.
The avoidance functional essentially contains all the information about $\eta_\lambda$ and, in fact, from Renyi's Theorem and its
generalizations (see, e.g., Theorem 9.2.XII in \cite{DaleyVereJOnes}) it is known that, with great generality, the distribution of
a simple point process is determined by the values of the mean avoidance function $A\mapsto P(\eta(A)=0)$ (we note that the term
\textit{avoidance function} is usually applied to this last functional; here we find it more convenient to make this change).
We think of $x+\lambda^{-1/d}Q(s)$ as a suitably scaled neighborhood of $x$ of a particular shape. We show in this paper 
that the weighted averages of avoidance functionals like the one in (\ref{functionaldefinition}) satisfy (see Theorems \ref{mainres}
and \ref{mainresK} below) 
\begin{equation}\label{mainbound0}
d_W\left(\frac{F_\lambda(\eta_\lambda)-\mbox{E}(F_\lambda(\eta_\lambda))}{\sqrt{\mbox{Var}(F_\lambda(\eta_\lambda))}},X \right)\leq
\frac 1{\sqrt{\lambda}} \frac{C_1}{C_2(\lambda)},\quad \lambda >0,
\end{equation}
where $X$ is a standard normal random variable, $C_1=C_1(\{Q(s)\}_s)$ is a constant that
measures the size of certain four-wise intersections of translates of the sets $B(s)$
and $C_2(\lambda)=C_2(\lambda;\{Q(s)\}_s)$ grows with $\lambda$ to a positive constant $C_2(\{Q(s)\}_s)$ which measures
measures the size of pairwise intersections of translates of the sets $B(s)$, with a similar result for $d_K$.
Thus, CLT's with rates for weighted averages of avoidance functionals can be proved from simple
geometrical considerations about the family $\{Q(s)\}$. 
While our approach lacks the generality of the results in \cite{LastPeccatiSchulte2014}, the particular
form of Mehler's formula for the class of functionals considered here allows for a more direct
treatment of the different terms that have to be controlled. Specifically, we do not use a general Poincaré
inequality but rely, instead, on an elementary bound for covariances of certain functionals or Poisson
random measures (see Lemma \ref{momentosPoisson} below). As a result we obtain bounds that are not covered 
by \cite{LastPeccatiSchulte2014} and offer a transparent geometric interpretation, see Remark \ref{comparison}
below for further details.

\medskip
We consider also the case in which we replace $\eta_\lambda$ by an empirical measure $\nu_n=\sum_{i=1}^n \delta_{X_i}$,
where $\{X_i\}$ are i.i.d. $\mathbb{R}^d$- valued r.v.'s. We show that, under suitable assumptions, for each
$\delta\in(0,\frac 1 4)$ there exists a positive constant, $C(\delta)$, such that
\begin{equation}\label{mainempir0}
d_W\left(\frac{F_n(\nu_n)-\mbox{E}(F_n(\nu_n))}{\sqrt{\mbox{Var}(F_n(\nu_n))}},X \right)\leq
\frac {C(\delta)}{n^{1/4-\delta}},\quad n\geq 1.
\end{equation}

\medskip
As an illustration of the power of (\ref{mainbound0}) and (\ref{mainempir0}) we show its application in
two classical problems in stochastic geometry: germ-grain models and quantization.
Germ-grain models are a common model in stochastic geometry (see, e.g., \cite{Stoyanetal}). 
They deal with random sets which arise as the union of sets centered around random
centers. A functional of interest in this setup is the volume of the resulting
random set. There CLT's for this random volume in the literature, see \cite{Penrose2007}
and the references therein. Here we consider random sets centered around points of
a Poisson process with unit intensity on $\mathbb{R}^d$, $\eta$, and
provide Berry-Essen bounds for the functional
\begin{equation}\label{germgrain0}
G_\lambda(\eta)=\ell\Big(\big(\cup_{z\in\eta} B(z,t)\big) \cap \big[-{\textstyle \frac 
{\lambda^{1/d}}2},{\textstyle \frac {\lambda^{1/d}}2}\big]^d \Big),
\end{equation}
where $B(z,t)$ is the ball of radius $t>0$ around $z$. We show (Theorem \ref{germTheo} in Section
4) that there is a constant $C(d,t)$ depending only on $t$ and $d$ such that
$$d\left(\frac{G_\lambda(\eta)-\mbox{E}(G_\lambda(\eta))}{\sqrt{\mbox{Var}(G_\lambda(\eta))}},X \right)\leq
\frac {C(d,t)}{\sqrt{\lambda}},\quad \lambda >0,$$
where $d=d_W$ or $d_K$.
We further provide a Berry-Esseen bound for the volume of the union of balls centered 
around the points of a empirical measure.

Our second application deals with the problem of quantization, that is, of 
approximation of a continuous measure by 
another measure concentrated on a finite set is a classical problem, originating 
in the information theory community, 
around the goal of signal discretization for an optimal signal transmission (see 
\cite{Zador1982}), which came to attract interest from other fields such as 
cluster analysis (see \cite{GrafLuschgy}) or finance (see, e.g. \cite{PagesWilbert2012}) among others. 
There are different approaches to the problem, depending on the way in which
the locations of the finitely supported probability are chosen. Here we focus on
quantization around random locations. More precisely, we consider the functional
\begin{equation}\label{quanterror0}
H(\eta_\lambda)=\int_{[0,1]^d} \min_{z\in\eta_\lambda}\|x-z \|^p dx=\sum_{z\in \eta_\lambda}
\int_{C(z,\eta_\lambda)}\|x-z \|^p dx,
\end{equation}
where $C(z,\eta_\lambda)$ is the Voronoi cell around $z$, that is, the set of points $x\in [0,1]^d$ which are
closer to $z$ than to any other point in the support of $\eta_\lambda$. A CLT for this and related functionals 
can be found in \cite{Yukich2008}. Here we prove (Theorem \ref{TeoremaPrincipal} below)
that there are positive constants $C,\lambda_0$ such that
$$d_{W}\left(\frac{H(\eta_\lambda)-\mbox{E}(H(\eta_\lambda))}{\sqrt{\mbox{Var}(H(\eta_\lambda))}},
X \right)\leq \frac{C}{\sqrt{\lambda}}, \lambda\geq \lambda_0,$$
where $X$ denotes a standard normal random variable.

The remaining sections of the paper are organized as follows. Section 2 provides a succint
summary of the main facts about the Malliavin calculus for Poisson functionals, providing
closed form expressions for the action of the key operators on avoidance functionals. Section
3 contains the main results of the paper, including the Berry-Esseen bounds for avoidance
functionals of Poisson process. The focus is put on the case of homogeneous Poisson processes, for
which the presentation is simpler, but we discuss the extension to general Poisson process and also to
empirical measures. Section 4 presents the application to the germ-grain model and to random quantization
errors. Finally, an Appendix provides some technical results used in some proofs in Section 3.

\section{Malliavin calculus for avoidance functionals}

As noted in the Introduction, the main result in this paper relies on the 
Malliavin calculus approach to Stein's method for Poisson functionals introduced
in \cite{Peccatietal2010}, which provides an upper bound for the Wasserstein distance
between the distribution of a Poisson functional and the standard normal distribution.
The upper bound is expressed in terms of certain Malliavin type operators, which, in turn,
are defined in terms of a Poisson chaos decomposition. For the sake of readibility we summarize
in this section the main facts about this chaos decomposition and the Malliavin type operators
and refer to \cite{LastPenrose2011} and \cite{Peccatietal2010} for further details. Later in the section
we particularize to the weighted averages of avoidance functionals, computing their chaos expansion
and the action of the relevant Malliavin operators on it.

\bigskip
We proceed then assuming that $\eta$ is a Poisson random measure on $\mathbb{R}^d$ with intensity measure $\mu$
(the theory is indeed much more general but this is enough for our purposes). A fundamental result is that any 
square integrable functional $F(\eta)$ admits a chaos decomposition, namely, it can be expressed as
\begin{equation}\label{PoissonChaos}
F(\eta)=\mbox{E}(F(\eta)) +\sum_{n=1}^\infty I_n(f_n)
\end{equation}
in $L_2$ sense, where $f_n(x_1,\ldots,x_n)$ is a real valued, square integrable function with respect to $\mu^n$,
symmetric in its $n$ arguments and $I_n$ denotes the multiple Wiener-It\^{o} integral with respect
to the compensated Poisson process, see, e.g. \cite{Peccatietal2010} or \cite{LastPenrose2011} for details.
For symmetric $f\in L_2((\mathbb{R}^d)^n, \mu^n)\cap L_1((\mathbb{R}^d)^n, \mu^n)$ we have
\begin{equation}\label{stochintegral}
I_n(f)=\sum_{I\subset [n]} (-1)^{n-|I|} \int_{(\mathbb{R}^d)^n} f(x_1,\ldots,x_n)d(\eta^{(I)}\otimes \mu^{([n]-I)})(x_1,\ldots,x_n)
\end{equation}
where $[n]$ is short for $\{1,\ldots,n\}$ and integration with respect to $\eta^{(I)}\otimes \mu^{([n]-I)}$ means summing $f(x_1,\ldots,x_n)$ 
with the $x_i$'s in the $I$ positions being all different and ranging in the support of $\eta$ and integrating with 
respect to $\mu$ all the other variables. We refer again to \cite{LastPenrose2011} for details. 

\medskip
The terms in decomposition (\ref{PoissonChaos}) are uncorrelated. In fact, for symmetric $f\in L_2((\mathbb{R}^d)^n,\mu^n)$,
$g \in L_2((\mathbb{R}^d)^m,\mu^m)$ we have
$$\mbox{E}(I_n(f)I_m(g))=I(m=n) n! \langle f,g \rangle_n,$$
where $\langle f,g \rangle_n$ denotes the usual inner product in $L_2((\mathbb{R}^d)^m,\mu^m)$. In particular, we see from (\ref{PoissonChaos})
that
\begin{equation}\label{formulavarianza}
\mbox{Var}(F(\eta))=\sum_{n=1}^\infty n! \|f_n\|_n^2
\end{equation}
and also that the functions $f_n$ in the chaos expansion (\ref{PoissonChaos}) are unique. In fact, these functions can be 
expressed in a remarkably simple way in terms of the difference operator, $D$, defined by the equation
\begin{equation}\label{differenceoperator}
D_{z}F(\eta)=F(\eta+\delta_{z})-F(\eta),
\end{equation}
with $\delta_z$ meaning Dirac's measure on $z$. Thus, $D_zF$ is the difference between the functional evaluated on the random set given 
by the support of $\eta$ supplemented by $z$ and its evaluation on the support of $\eta$. In terms of this operator,
the functions $f_n$ in (\ref{PoissonChaos}) are (up to a scaling constant) the expected values of the iterated difference
operator acting on $F(\eta)$. More precisely,
\begin{equation}\label{formulaparafn}
f_n(x_1,\ldots,x_n)=\frac 1 {n!}\mbox{E}(D_{x_1,\ldots,x_n}F(\eta)).
\end{equation}
Here $D_{x_1,\ldots,x_n}F(\eta)=D_{x_1} D_{x_2,\ldots,x_n}F(\eta)$.
A useful fact about the operator $D_{x_1,\ldots,x_n}$ is that
\begin{equation}\label{expressionDxn}
D_{x_1,\ldots,x_n}F(\eta)=\sum_{I\subset \{1,\ldots,n\}} (-1)^{n+|I|} F(\eta +\sum_{i\in I} \delta_{x_i}),
\end{equation}
see, e.g., \cite{Schulte2012}.

\medskip
Alternatively, the difference operator can be expressed in terms of an orthogonal expansion. In fact, for $F$ as 
in (\ref{PoissonChaos})
\begin{equation}\label{DzAlternative}
D_z F(\eta)= \sum_{n=1}^\infty n I_{n-1}(f_n(z,\cdot)).
\end{equation}
Other important Malliavin operators are expressed through similar orthogonal expansions. These include the
Ornstein-Uhlenbeck operator, $L$, given by
$$L F(\eta)=-\sum_{n=1}^\infty n I_{n}(f_n),$$
and the inverse Ornstein-Uhlenbeck operator, $L^{-1}$, for which
$$L^{-1} F(\eta)=-\sum_{n=1}^\infty \frac 1 n I_{n}(f_n).$$
The domain of $L$ consists of the $L_2$ functionals such that $\mbox{Var}(F(\eta))=\sum_{n=1}^\infty  n\, n! \|f_n\|_n^2<\infty$
while that of $L^{-1}$ is the set of centered $L_2$ functionals. 

\medskip
The operator $D_z L^{-1}$ will be of special interest for us. For $F$ in the domain of $L^{-1}$ we have
\begin{equation}\label{otramas}
-D_z L^{-1} F=\sum_{n=1}^\infty I_{n-1} (f_n(z,\cdot)).
\end{equation}
Our next results provide a simple closed form for $D_z L^{-1}F(\eta)$ in the case
of avoidance functionals. Part d) in the next Lemma can be obtained from 
Mehler's formula (Theorem 3.2 in\cite{LastPeccatiSchulte2014}). Here we
provide a self-contained elementary proof.

\begin{Lemma}\label{L1Indicator}
If $\eta$ is a Poisson point process on $\mathbb{R}^d$ with intensity measure $\mu$,
$A\subset \mathbb{R}^d$ a Borel set such that $\mu(A)<\infty$ and $F(\eta)=\mathbf{1}(\eta(A)=0)-e^{-\mu(A)}$ then,
\begin{itemize}
\item[(a)] the $f_n$ functions in the chaos expansion (\ref{PoissonChaos}) are given by
\begin{equation}\label{formulafn1}
f_n(x_1,\ldots,x_n)=\frac {(-1)^n} {n!} e^{-\mu(A)} \prod_{i=1}^n \mathbf{1}_A(x_i),\quad n\geq 1.
\end{equation}
\item[(b)] $$I_n(f_n)=\frac{e^{-\mu(A)}}{n!}\sum_{k=0}^n (-1)^k \binom{n}{k}\binom{\eta(A)}{k} k! \mu(A)^{n-k} ,$$
\item[(c)] $$D_z F(\eta)=-\mathbf{1}_A(z) F(\eta)\qquad \mbox{and}$$
\item[(d)] $$D_z L^{-1}F(\eta)=\mathbf{1}_A(z) \int_0^1 t^{\eta(A)} e^{-\mu(A)t}dt.$$
\end{itemize}
\end{Lemma}

\medskip
\noindent 
\textbf{Proof.} We compute first the $f_n$ kernels in the chaos expansion (\ref{PoissonChaos}). 
For convenience, given a point measure $\rho$ we rewrite the set $(\rho(A)=0)$ as 
$(A\cap \rho =\emptyset)$. From (\ref{expressionDxn}) we see that
$$f_n(x_1,\ldots,x_n)=\frac 1 {n!}\sum_{I\subset \{1,\ldots,n\}} (-1)^{n+|I|} 
\mbox{P}\left(A\cap (\eta \cup (\cup_{i\in I} \{x_i\}))=\emptyset\right).$$
Now, if $x_i\in A$ for all $i\in\{1,\ldots,n\}$ then all the
terms in the expansion of $f_n(x_1,\ldots,x_n)$ vanish except for the case $I=\emptyset$ and,
consequently, $f_n(x_1,\ldots,x_n)=\frac{(-1)^n} {n!} \mbox{P}\left(A\cap \eta=\emptyset\right)=
\mbox{P}\left(\eta(A)=\emptyset\right)=
\frac{(-1)^n} {n!} e^{-\mu(A)}.$ Assume, on the contrary, that there is some index $i\in\{1,\ldots,n\}$
such that $x_i\notin A$. Then, for every $J\subset\{1,\ldots,n\}\backslash \{i\}$ we have
$A\cap (\eta \cup (\cup_{j\in J} \{x_j\}))= A\cap (\eta \cup (\cup_{j\in J\cup \{i \}} \{x_j\}))$. 
Hence,
\begin{eqnarray*} 
\lefteqn{f_n(x_1,\ldots,x_n)=\frac 1{n!}\sum_{J\subset \{1,\ldots,n\}\backslash \{ i\} } (-1)^{n+|J|} \left[
\mbox{P}\left(A\cap (\eta \cup (\cup_{j\in J} \{x_j\}))=\emptyset\right) \right.}\hspace*{6.5cm} \\
&& \left. - \mbox{P}\left(A\cap (\eta \cup (\cup_{j\in J\cup \{i \}} \{x_j\}))=\emptyset\right)\right]= 0.
\end{eqnarray*}
Combining these facts we obtain (\ref{formulafn1}).

\medskip
We turn now to (b). Fix $I\subset [n]$ with $|I|=k$, say. Then
$$\int_{(\mathbb{R}^d)^n} f(x_1,\ldots,x_n)d(\eta^{(I)}\otimes \mu^{([n]-I)})(x_1,\ldots,x_n)=\frac {(-1)^n} {n!} e^{-\mu(A)} \mu(A)^{n-k} k! 
 \binom{\eta(A)}{k}$$
and (\ref{stochintegral}) becomes
$$I_n(f_n)=\frac{e^{-\mu(A)}}{n!}\sum_{k=0}^n (-1)^k \binom{n}{k} k! \mu(A)^{n-k}  \binom{\eta(A)}{k},$$
proving (b). Part (c) is obvious. Finally, to prove (d) using (\ref{otramas}) and linearity of the stochastic integral we see that
\begin{eqnarray*}
D_z L^{-1}F(\eta)&=& \mathbf{1}_A(z)e^{-\mu(A)} \sum_{n=1}^{\infty} \frac{(-1)^{n-1}}{n!} \left[ \sum_{k=0}^{n-1} \binom{n-1}{k} (-1)^{n-1-k}  \mu(A)^{n-1-k} k!
\binom{\eta(A)}k\right]\\
&=& \mathbf{1}_A(z)e^{-\mu(A)} \sum_{m=0}^{\infty} \frac{(-1)^{m}}{(m+1)!} \left[ \sum_{k=0}^{m} \binom{m}{k} (-1)^{m-k}  \mu(A)^{m-k} k!
\binom{\eta(A)}k\right]\\
&=& \mathbf{1}_A(z)e^{-\mu(A)} \sum_{k=0}^{\infty} (-1)^k \binom{\eta(A)}k \left[ \sum_{m=k}^{\infty} \frac{1}{m+1}  \frac{\mu(A)^{m-k}}{(m-k)!}
\right]\\
&=& \mathbf{1}_A(z)e^{-\mu(A)} \sum_{k=0}^{\infty} (-1)^k \binom{\eta(A)}k \left[ \sum_{n=0}^{\infty} \frac{1}{n+k+1}  \frac{\mu(A)^{n}}{(n)!}
\right]\\
&=& \mathbf{1}_A(z)e^{-\mu(A)} \int_0^1 \left[\sum_{k=0}^{\infty}  \binom{\eta(A)}k (-x)^k \left[ \sum_{n=0}^{\infty}   \frac{(x\mu(A))^{n}}{(n)!}
\right]\right]dx\\
&=& \mathbf{1}_A(z)e^{-\mu(A)} \int_0^1 \left[\sum_{k=0}^{\infty}  \binom{\eta(A)}k (-x)^k \right] e^{\mu(A) x} dx\\
&=& \mathbf{1}_A(z) \int_0^1 (1-x)^{\eta(A)} e^{-\mu(A)(1- x)} dx.
\end{eqnarray*}
Here we have used that $\frac{1}{n+k+1}=\int_{0}^1 x^{k+n}dx$. The fact that $\mbox{E}((1-x)^{\eta(A)})=e^{-\mu(A)x}$ shows that the last 
integral is a.s. finite and justifies the exchange of integration/sums in the above lines. The change of variable $t=1-x$ completes the proof.
\quad $\Box$

\medskip 
From this Lemma \ref{L1Indicator} we obtain a simple description of the action of the
operators $D$ and $DL^{-1}$ over the avoidance functionals defined in (\ref{functionaldefinition}). 
The proof is a simple consequence of Lemma \ref{L1Indicator} and linearity of the 
operators $D$ and $DL^{-1}$. We omit details.

\begin{Corollary}\label{actionoperators}
Assume $\eta$ is a Poisson point process on $\mathbb{R}^d$ with intensity measure $\mu$
defined on the probability space $(\Omega,\mathcal{F},P)$, $\mathcal{A}$ a bounded 
open set on $\mathbb{R}^d$ and $(\mathcal{B},\mathcal{G},\nu)$
a measure space. Assume further that, for each $s\in\mathcal{B}$, $Q(s)$ is a Borel set on $\mathbb{R}^d$
chosen in such a way that the map $(\omega,x,s)\mapsto \mathbf{1}(\eta (x+Q(s))=0)$
is $\mathcal{F}\otimes \beta^d \otimes \beta$ measurable. Then, if 
$$
F(\eta)=\int_{\mathcal{A}\times\mathcal{B}} \mathbf{1}(\eta (x+Q(s))=0) d(\ell\otimes \nu)(x,s)$$
and 
$$\int_{(\mathcal{A}\times\mathcal{B})^2} e^{-\mu(Q(x,s)\cup Q(y,t))} dxd\nu(s)dyd\nu(t)<\infty, $$
we have 
$$D_zF(\eta)=-\int_{\mathcal{A}\times\mathcal{B}} \mathbf{1}(\eta (x+Q(s))=0)\mathbf{1}(z\in x+Q(s)) d(\ell\otimes \nu)(x,s),$$
and
$$D_zL^{-1}(F(\eta)-E(F(\eta)))=\int_{\mathcal{A}\times\mathcal{B}} \Big(\int_0^1 t^{\eta (x+Q(s))}
e^{-\mu(x+Q(s))t}dt\Big)\mathbf{1}(z\in x+Q(s)) d(\ell\otimes \nu)(x,s).$$
\end{Corollary}

\section{Normal approximation for avoidance functionals}

In this section we provide an upper bound for the Wasserstein and Kolmogorov distance between the distribution of
the standardized avoidance functional and the standard normal distribution. Our approach is 
based on the following result, which is a simplified version of Theorem 3.1 in \cite{Peccatietal2010}
and Theorem 3.1 in \cite{EiselbacherThale2013}.

\begin{Theorem}\label{WassersteinSteinPoisson}
If $\eta$ is a Poisson point process with nonatomic, $\sigma$-finite 
intensity measure $\mu$, $F=F(\eta)$ satisfies $0<\mbox{\em Var}(F)<\infty$ and
$X\sim N(0,1)$ then
$$d_W\left(\frac{F-\mbox{\em E}(F)}{\sqrt{\mbox{\em Var}(F)}},X \right)\leq 
\frac{(\mbox{\em Var} (\int A(x)B(x) d\mu(x)))^{1/2}}{\mbox{\em Var}(F)}
+\frac{\mbox{\em E} (\int A(x)^2 |B(x)|  d\mu(x))}{(\mbox{\em Var}(F))^{3/2}}
$$
\begin{eqnarray*}
d_K\left(\frac{F-\mbox{\em E}(F)}{\sqrt{\mbox{\em Var}(F)}},X \right)&\leq &
\frac{(\mbox{\em Var} (\int A(x)B(x) d\mu(x)))^{1/2}}{\mbox{\em Var}(F)}+
{\textstyle \frac {4+\sqrt{2\pi}}{8}}\frac{\mbox{\em E} (\int A(x)^2 |B(x)|  d\mu(x))}{(\mbox{\em Var}(F))^{3/2}}\\
&+& {\textstyle \frac 1 2}\frac{(\mbox{\em Var}(\int A(x)^2 |B(x)|  d\mu(x) ))^{1/2}}{(\mbox{\em Var}(F))^{3/2}}\\
&+& \frac{(\mbox{\em E}(\int A(x)^2 B(x)^2  d\mu(x) +\int\int |D_y C(x) D_x C(y)|d\mu(x)d\mu(y))^{1/2}}{\mbox{\em Var}(F)}
\end{eqnarray*}
where $A(x)= D_x F$, $B(x)=-D_x L^{-1}(F-\mbox{\em E}(F))$ and $C(x)=A(x)|B(x)|$.
\end{Theorem}

\medskip
\noindent \textbf{Proof.} The first inequality is simply Theorem 3.1 in \cite{Peccatietal2010} applied to 
$\tilde{F}=\frac{F-\mbox{\scriptsize E}(F)}{\sqrt{\mbox{\scriptsize Var}(F)}}$,
after noting that $\mbox{E}(\int A(x)B(x) d\mu(x))=\mbox{Var} (F)$. For the second we use Theorem 3.1 in  
\cite{EiselbacherThale2013} to get
\begin{eqnarray*}
d_K\left(\frac{F-\mbox{E}(F)}{\sqrt{\mbox{Var}(F)}},X \right)&\leq &
\frac{(\mbox{Var} (\int A(x)B(x) d\mu(x)))^{1/2}}{\mbox{Var}(F)}+
{\textstyle \frac {\sqrt{2\pi}}{8}}\frac{\mbox{E} (\int A(x)^2 |B(x)|  d\mu(x))}{(\mbox{Var}(F))^{3/2}}\\
&+&{\textstyle \frac {1}{2}}\frac{\mbox{E} \big(|F-\mbox{E}(F)|\int A(x)^2 |B(x)|  d\mu(x)\big)}{(\mbox{Var}(F))^{2}}\\
&+& \frac{\sup_{t\in\mathbb{R}}\mbox{E}\big(\int A(x) |B(x)| D_x\mathbf{1}(F>t) d\mu(x) \big)}{\mbox{Var}(F)}.
\end{eqnarray*}
For the third term in the last upper bound we use Schwarz's inequality and the fact that $\sqrt{a+b}\leq \sqrt{a}+\sqrt{b}$ for $a,b\geq0$,
to obtain
\begin{eqnarray}\label{dK1}
\lefteqn{\frac{\mbox{E} \big(|F-\mbox{E}(F)|\int A(x)^2 |B(x)|  d\mu(x)\big)}{(\mbox{Var}(F))^{2}}\leq 
\frac{\Big(\mbox{E} \big(\int A(x)^2 |B(x)|  d\mu(x)\big)^2\Big)^{1/2}}{(\mbox{Var}(F))^{3/2}}}\hspace*{3cm}\\
&=&\nonumber
\frac{\Big(\mbox{Var} \big(\int A(x)^2 |B(x)|  d\mu(x)\big)+ \big(
\mbox{E} \int A(x)^2 |B(x)|  d\mu(x)\big)^2\Big)^{1/2}}{(\mbox{Var}(F))^{3/2}}\\
&\leq & \nonumber 
\frac{\big(\mbox{Var} \big(\int A(x)^2 |B(x)|  d\mu(x)\big)\big)^{1/2}+ 
\mbox{E} \int A(x)^2 |B(x)|  d\mu(x)}{(\mbox{Var}(F))^{3/2}}.
\end{eqnarray}
To deal with the fourth term we assume that the corresponding term in the statement
is finite (there is nothing to prove otherwise). Then we can use Lemma 2.2 and Proposition 2.3 
in \cite{LastPeccatiSchulte2014} and Schwarz's inequality to see that $C$ is in $\mbox{dom}(\delta)$,
with $\delta$ the Skorohod integral operator, and
\begin{eqnarray}\nonumber
\lefteqn{\mbox{E}\big({\textstyle\int } A(x) |B(x)| D_x\mathbf{1}(F>t) d\mu(x) \big)=\mbox{E}(\mathbf{1}(F>t) \delta(C))\leq \big( 
\mbox{E}(\mathbf{1}(F>t) \delta(C)^2)\big)^{1/2}}\hspace*{1cm}\\
&\leq &
\mbox{E}(\delta(C)^2)\big)^{1/2}\leq 
(\mbox{\em E}({\textstyle\int} A(x)^2 B(x)^2  d\mu(x) +{\textstyle\int\int} |D_y C(x) D_x C(y)|d\mu(x)d\mu(y))^{1/2}.
\end{eqnarray}
This completes the proof. \hfill $\Box$

\medskip
We plan to apply now Theorem \ref{WassersteinSteinPoisson} to avoidance functionals.
The next result will be essential for this task.

\begin{Lemma}\label{momentosPoisson}
If $\eta$ is a Poisson random measure with intensity measure $\mu$ and $A,B,C$ and $D$ are
measurable sets of  finite $\mu$-measure then
$$\mbox{\em E}\left(\mathbf{1}(\eta(A)=0) \int_0^1 u^{\eta(B)}e^{-\mu(B)u}du \right)= e^{-\mu(A\cup B)} 
\frac{1-e^{-\mu(A\cap B)}}{\mu(A\cap B)},$$
($=e^{-\mu(A\cup B)}$ if $\mu(A\cap B)=0$) and
\begin{eqnarray*}
\lefteqn{\mbox{\em Cov}\left(\mathbf{1}(\eta(A)=0) \int_0^1 u^{\eta(B)}e^{-\mu(B)u}du,
\mathbf{1}(\eta(C)=0) \int_0^1 v^{\eta(D)}e^{-\mu(D)v}dv \right)}\hspace*{4cm}\\
&&\leq e^{-\mu(A\cup B \cup C \cup D)} 
\mathbf{1}((A\cup B)\cap (C\cup D) \ne \emptyset).
\end{eqnarray*}
\end{Lemma}

\medskip
\noindent 
\textbf{Proof.} From the fact $\mathbf{1}(\eta(A)=0) u^{\eta(B)}=\mathbf{1}(\eta(A)=0) u^{\eta(B\cap A^C)}$
and independence of $\eta (A)$ and $\eta(B\cup A^C)$
we see that
\begin{eqnarray*}
\mbox{\em E}\left(\mathbf{1}(\eta(A)=0) \int_0^1 u^{\eta(B)}e^{-\mu(B)u}du \right)&=&e^{-\mu(A)}
\int_0^1 \mbox{E}(u^{\eta(B\cap A^C)})e^{-\mu(B)u}du\\
&=&e^{-\mu(A)}
\int_0^1 e^{-\mu(B\cap A^C)(1-u)}e^{-\mu(B)u}du\\
&=&e^{-\mu(A\cup B)}
\int_0^1 e^{-\mu(A\cap B)u}du
\end{eqnarray*}
and prove the first part of the Lemma. For the upper bound for the covariance
we note that we have independence (hence null covariance) if $(A\cup B)\cap (C\cup D)=\emptyset$. On the other hand,
arguing as above we see that
\begin{eqnarray*}
\lefteqn{\mbox{E}\left(\mathbf{1}(\eta(A)=0) \int_0^1 u^{\eta(B)}e^{-\mu(B)u}du
\mathbf{1}(\eta(C)=0) \int_0^1 v^{\eta(D)}e^{-\mu(D)v}dv \right)}\hspace*{0cm}\\
&&=\mbox{E}\left(\mathbf{1}(\eta(A\cup C)=0) \int_{(0,1)^2} u^{\eta(B)}v^{\eta(D)}e^{-\mu(B)u}e^{-\mu(D)v}du dv\right)\hspace*{5.3cm}
\end{eqnarray*}
\begin{eqnarray*}
&&=e^{-\mu(A\cup C)} \int_{(0,1)^2}\mbox{E}( u^{\eta( A^C\cap B \cap C^C \cap D^C)} v^{\eta( A^C\cap B^C\cap C^C \cap D)}
(uv)^{\eta( A^C\cap B\cap C^C \cap D)})
e^{-\mu(B)u}e^{-\mu(D)v}du dv\\
&&=e^{-\mu(A\cup C)} \int_{(0,1)^2}e^{-\mu( A^C\cap B \cap C^C \cap D^C)(1-u)}e^{-\mu( A^C\cap B^C\cap C^C \cap D)(1-v)}
e^{-\mu( A^C\cap B\cap C^C \cap D)(1-uv) }\\
&& \hspace*{7cm}\times
e^{-\mu(B)u}e^{-\mu(D)v}du dv\\
&&=e^{-\mu(A \cup B \cup C \cup D)}\int_{(0,1)^2} e^{-\mu( B \cap (A\cup C\cup D)  )u-\mu( D \cap (A\cup B\cup C))v
+\mu( A^C\cap B\cap C^C \cap D)uv}dudv\\
&&\leq e^{-\mu(A \cup B \cup C \cup D)},
\end{eqnarray*}
where the last upper bound follows from the fact that the exponent in the last integral is nonpositive for every $u,v\in(0,1)$.\quad $\Box$

\medskip
Now we focus on the following setup. $\eta_\lambda$ will denote a homogenous Poisson process on $\mathbb{R}^d$
with intensity $\lambda$ defined on $(\Omega,\mathcal{F},P)$, $\mathcal{A}$ a nonempty, bounded open set on $\mathbb{R}^d$ and $(\mathcal{B},\mathcal{G},\nu)$
a measure space and will assume that, for each $s\in\mathcal{B}$, $Q(s)$ is a Borel set on $\mathbb{R}^d$
chosen in such a way that the map $(\omega,x,s)\mapsto \mathbf{1}(\eta_\lambda (x+\lambda^{-1/d}Q(s))=0)$
is $\mathcal{F}\otimes \beta^d \otimes \mathcal{G}$ measurable and such that
\begin{equation}\label{fintemean}
\int_{\mathcal{B}}  e^{-\ell(Q(s))}d\nu(s)<\infty.
\end{equation}
We will write $Q_\lambda(x,s)=x+\lambda^{-1/d}Q(s)$ 
and $Q(x,s)=Q_1(x,s)$.
Now, we consider the avoidance functional
\begin{equation}\label{functionaldefinition2}
F_\lambda(\eta_\lambda)=\int_{\mathcal{A}\times\mathcal{B}}\mathbf{1}(\eta_\lambda (Q_\lambda(x,s))=0) d(\ell\otimes \nu)(x,s).
\end{equation}
Our next results provides simple expressions for the mean and variance of $F_\lambda(\eta_\lambda)$.
For a simpler statement we introduce the functions
$$V(z):=\int_{\mathcal{B}\times\mathcal{B}}
e^{-(\ell(Q(0,s))+\ell(Q(z,t))) }\big[ e^{\ell (Q(0,s)\cap (Q(z,t)))}-1 \big]d\nu(s)d\nu(t), \quad z\in\mathbb{R}^d,$$
and 
\begin{equation}\label{C1lambda}
C_1(\lambda):=\int_{U_\lambda} V(z)dxdz,\quad \lambda>0,
\end{equation}
where $U_\lambda=\{(x,z): \, x\in \mathcal{A},z\in \lambda^{1/d}(\mathcal{A}-x))  \}$. We observe that $C_1(\lambda)$
is a nondecreasing function that satisfies 
\begin{equation}\label{C1}
\lim_{\lambda \to\infty} C_1(\lambda)=C_1:=\ell(\mathcal{A}) \int_{\mathbb{R}^d} V(z) dz.
\end{equation}

\begin{Lemma}\label{moments}
Under (\ref{fintemean}), if $F_\lambda(\eta_\lambda)$  is the 
avoidance functional defined in (\ref{functionaldefinition2}) and the constant $C_1$ in
(\ref{C1}) is finite then $F_\lambda(\eta_\lambda)$ has finite second moment and
$$\mbox{\em E}(F_\lambda(\eta_\lambda))=\ell(\mathcal{A})\int_{\mathcal{B}}  e^{-\ell(Q(s))}d\nu(s)$$
and 
$$\lambda \mbox{\em Var}(F_\lambda(\eta_\lambda))=C_1(\lambda),$$
with $C_1(\lambda)$ defined by (\ref{C1lambda}).
\end{Lemma}

\medskip
\noindent 
\textbf{Proof.} 
For computing the mean we use Fubini's theorem and the fact that $\lambda \ell(Q_\lambda(x,s))=\ell(Q(s))$ to get
\begin{eqnarray*}
\mbox{E}(F(\eta_\lambda))
&=& \int_{\mathcal{A}\times\mathcal{B}}  \mbox{P}(\eta_\lambda(Q_\lambda(x,s))=0 ) dxd\nu(s)\\
&=& \int_{\mathcal{A}\times\mathcal{B}} e^{-\lambda\ell(Q_\lambda(x,s))}   dxd\nu(s)=\ell(\mathcal{A})\int_{\mathcal{B}}  e^{-\ell(Q(s))}d\nu(s).
\end{eqnarray*}
Turning now to the
variance we use again Fubini's theorem to obtain
\begin{eqnarray}\nonumber
\mbox{E}(F(\eta_\lambda)^2)&=&\int_{\mathcal{A}\times\mathcal{A} }\left[ \int_{\mathcal{B}\times\mathcal{B}} 
\mbox{P}(\eta_\lambda(Q_\lambda(x,s)\cup Q_\lambda(y,t))=0 ) d\nu(s)d\nu(t)\right]dx dy\\
\nonumber
&=&\int_{\mathcal{A}\times\mathcal{A} }\left[ \int_{{\mathcal{B}\times\mathcal{B} }} 
e^{-\lambda \ell(Q_\lambda(x,s)\cup Q_\lambda(y,t))} d\nu(s)d\nu(t)\right]dx dy\\
\label{momento2} 
&=&\int_{\mathcal{A}\times\mathcal{A} }\left[ \int_{{\mathcal{B}\times\mathcal{B} }} 
e^{-\lambda \ell(Q_\lambda(0,s)\cup Q_\lambda(y-x,t))} d\nu(s)d\nu(t)\right]dx dy,
\end{eqnarray}
the last equality coming from translation invariance of Lebesgue measure. The change of variable 
$y=x+\lambda^{-\frac 1 d} z$ and the fact that $\lambda\ell(Q_\lambda(0,s)\cup Q_\lambda(\lambda^{-1/d}z,t))=
\ell(Q(0,s)\cup Q(z,t))$ yield
\begin{eqnarray*}
\mbox{E}(F(\eta_\lambda)^2)=\lambda^{-1} \int_{U_\lambda } \left[\int_{\mathcal{B}\times\mathcal{B} }
 e^{-\ell(Q(0,s)\cup Q(z,t))}d\nu(u)\,d\nu(v)    \right]dx dz.
\end{eqnarray*}
A similar computation shows that
$$(\mbox{E}(F(\eta_\lambda)))^2 
=\lambda^{-1} \int_{U_\lambda } \left[\int_{\mathcal{B}\times\mathcal{B} }
 e^{-\ell(Q(0,s))+\ell(Q(z,t))}d\nu(u)\,d\nu(v)    \right]dx dz.$$
and, combining the last two equations, we conclude that
$$\mbox{Var}(F(\eta_\lambda))=\lambda^{-1}\int_{U_\lambda } V(z)dx dz.$$
\quad $\Box$

\medskip

We are ready now for the main results of this section, namely, explicit upper bounds for the 
Wasserstein or Kolmogorov distance between the law of standardized Poisson avoidance functionals and the
standard normal distribution. We  consider first the Wasserstein case and we introduce the constants
\begin{eqnarray}\nonumber
C_{2,a}&=&\ell(\mathcal{A})\int_{(\mathbb{R}^d)^2} W_a(z_1,z_2) dz_1 dz_2,\\
C_{2,b}&=&\ell(\mathcal{A})\int_{(\mathbb{R}^d)^3} W_b(z_1,z_2,z_3) dz_1 dz_2 dz_3\label{constantsC2}
\end{eqnarray}
with
$$W_a(z_1,z_2)=\int_{\mathcal{B}^3}  e^{-\ell(\cup_{i=0}^2 Q_i)}  \ell(\cap_{i=0}^2Q_i) \prod_{i=0}^2 \nu(ds_i),$$
$$W_b(z_1,z_2,z_3)=\int_{\mathcal{B}^4}  e^{-\ell(\cup_{i=0}^3 Q_i)} \ell(Q_0\cap Q_1) \ell(Q_2\cap Q_3) 
\mathbf{1}([B_0\cup B_1]\cap [B_2\cup B_3]\ne \emptyset) \prod_{i=0}^3\nu(ds_i),$$
where, in the last two integrals $Q_i$ is short notation for $Q(z_i,s_i)$ and $z_0=0$. Note that
$C_{2,a}$ and $C_{2,b}$ weighted averages of the size of different three- or four-wise intersections of translates of the sets $Q(s)$.
As we show next, control of these averages is all that is needed to give CLT's with rates in Wasserstein distance for 
avoidance functionals.
\begin{Theorem}\label{mainres}
Assume $\eta_\lambda$ is a homogenous Poisson process on $\mathbb{R}^d$
with intensity $\lambda$ defined on $(\Omega,\mathcal{F},P)$, $\mathcal{A}$ a nonempty, bounded open set on $\mathbb{R}^d$ and $(\mathcal{B},\mathcal{G},\nu)$
a measure space. Assume further that, for each $s\in\mathcal{B}$, $Q(s)$ is a Borel set on $\mathbb{R}^d$
chosen in such a way that the map $(\omega,x,s)\mapsto \mathbf{1}(\eta_\lambda (x+\lambda^{-1/d}Q(s))=0)$
is $\mathcal{F}\otimes \beta^d \otimes \mathcal{G}$ measurable. If 
$$F_\lambda(\eta)=\int_{\mathcal{A}\times\mathcal{B}} \mathbf{1}(\eta (x+\lambda^{-1/d}Q(s))=0) d(\ell\otimes \nu)(x,s),$$
(\ref{fintemean}) holds and the constant $C_1$ in (\ref{C1}) is finite, then $F_\lambda(\eta_\lambda)$
has finite second moment and
$$d_W\left(\frac{F_\lambda(\eta_\lambda)-\mbox{\em E}(F_\lambda(\eta_\lambda))}{\sqrt{\mbox{\em Var}(F_\lambda(\eta_\lambda))}},X \right)\leq
\frac 1{\sqrt{\lambda}} \Big(\frac{C_{2,a}}{C_1(\lambda)^{3/2}}+\frac{C_{2,b}^{1/2}}{C_1(\lambda)}\Big),\quad \lambda >0,$$
where $X$ denotes a standard normal random variable and $C_{2,a}$ and $C_{2,b}$ are given by (\ref{constantsC2}).
\end{Theorem}

\medskip
\noindent\textbf{Proof.}
From Corollary \ref{actionoperators} and Fubini's theorem we see that
\begin{eqnarray*}
\lefteqn{\int_{\mathbb{R}^d} (D_z F_\lambda(\eta_\lambda))^2 |D_zL^{-1}(F-E(F)|) dz}\hspace*{1cm}\\
&=&\int_{(\mathcal{A}\times\mathcal{B})^3 } 
\mathbf{1}\big(\eta_\lambda\big( \cup_{i=0}^1 Q_\lambda(x_i,s_i)\big)=0\big)\big(
{\textstyle\int_0^1} u^{\eta_\lambda(Q_\lambda(x_2,s_2))} e^{-u\lambda \ell (Q_\lambda(x_2,s_2))}du\big)\\
&\times & \ell \big(
\cap_{i=0}^2 Q_\lambda(x_i,s_i)\big)d\prod_{i=0}^2(\ell\times \nu)(x_i,s_i)
\end{eqnarray*}
and, therefore, using Lemma \ref{momentosPoisson} (and the fact that $1-e^{-x}\leq x, x\geq 0$),
\begin{eqnarray*}
\lefteqn{E\big({\textstyle \int_{\mathbb{R}^d}} (D_z F_\lambda(\eta_\lambda))^2 |D_zL^{-1}(F_\lambda(\eta_\lambda)-E(F_\lambda(\eta_\lambda))| dz\big)}\hspace*{1cm}
\\
&\leq &\int_{(\mathcal{A}\times\mathcal{B})^3 } 
e^{-\lambda \ell\big( \cup_{i=0}^2 Q_\lambda(x_i,s_i)\big)} \lambda \ell \big(
\cap_{i=0}^2 Q_\lambda(x_i,s_i)\big)d\prod_{i=0}^2(\ell\times \nu)(x_i,s_i)\\
&=&\int_{\mathcal{A}^3} \Big[\int_{\mathcal{B}^3}
e^{-\lambda \ell\big( \cup_{i=0}^2 Q_\lambda(x_i,s_i)\big)} \lambda \ell \big(
\cap_{i=0}^2 Q_\lambda(x_i,s_i)\big)\prod_{i=0}^2d\nu(s_i) \Big]\prod_{i=0}^2 dx_i.
\end{eqnarray*}
Hence, if we change variables, $x_i=x_0+\lambda^{-1/d}z_i$, $i=1,2$, denote 
$V_\lambda=\{(x_0,z_1,z_2): \, x_0\in \mathcal{A},z_i\in \lambda^{1/d}(\mathcal{A}-x_0)), i=1,2  \}$,
and observe that $\lambda \ell\big( \cup_{i=0}^2Q_\lambda(x_i,s_i)\big)=\ell\big( \cup_{i=0}^2 Q(z_i,s_i)\big)$
(with $z_0=0$) and similarly for $\lambda \ell\big( \cap_{i=0}^2 Q_\lambda(x_i,s_i)\big)$,  we see that
\begin{equation}\label{finalbounda}
E\big({\textstyle \int_{\mathbb{R}^d}} (D_z F_\lambda(\eta_\lambda))^2 |D_zL^{-1}(F-E(F)| dz\big)=\frac 1{\lambda^2}
\int_{V_\lambda}W_a(z_1,z_2)dx_0dz_1dz_2\leq \frac {C_{2,a}}{\lambda^2}.
\end{equation}
In a similar fashion, we obtain from Corollary \ref{actionoperators} that
$\mbox{Var}(\langle DF_\lambda(\eta_\lambda), DL^{-1}(F_\lambda(\eta_\lambda)-\mbox{E}F_\lambda(\eta_\lambda)\rangle)$
equals the variance of 
$$\int_{\mathcal{A}^2\times\mathcal{B}^2} \mathbf{1}(\eta_\lambda(Q_\lambda(x_0,s_0))=0)\Big(\int_0^1 t^{\eta_\lambda(Q_\lambda(x_1,s_1))} 
e^{-\lambda \eta_\lambda(Q_\lambda(x_0,s_0))t}dt  \Big)\lambda\ell(\cap_{i=0}^1 Q_\lambda(x_i,s_i))\prod_{i=0}^1(dx_i  d\nu(s_i)),$$
which, using the covariance inequality of Lemma \ref{momentosPoisson}, is upper bounded by
$$
\int_{\mathcal{A}^4\times\mathcal{B}^4} e^{-\lambda \ell(\cup_{i=0}^3 Q_i)}\lambda\ell(\cap_{i=0}^1 Q_i)
\lambda\ell(\cap_{i=2}^3 Q_i)\mathbf{1}((\cup_{i=0}^1 Q_i)\cap (\cup_{i=2}^3 Q_i)\ne \emptyset)\prod_{i=0}^3(dx_i  d\nu(s_i)),$$
with $Q_i=Q_\lambda(x_i,s_i)$, $i=0,\ldots,3$ in this last integral. A change of variable as above yields now
\begin{equation}\label{finalboundb}
\mbox{Var}(\langle DF_\lambda(\eta_\lambda), DL^{-1}(F_\lambda(\eta_\lambda)-\mbox{E}F_\lambda(\eta_\lambda)\rangle)\leq
\frac {C_{2,b}}{\lambda^3}.
\end{equation}
Now, the conclusion follows from Theorem \ref{WassersteinSteinPoisson}, Lemma \ref{momentosPoisson} and
(\ref{finalbounda}).\quad $\Box$

\begin{Remark}\label{comparison}{\em
Theorem \ref{mainres} should be compared to Proposition 1.3 in \cite{LastPeccatiSchulte2014}.
We observe that we do not need, for instance, to assume uniform upper bounds for (moments of)
$D_zF_\lambda(\eta_\lambda)$ in order to have a CLT with rate $\lambda^{-1/2}$ in
Wasserstein distance (of course, the gain comes from the fact that we
have restricted ourselves to a particular class of functionals and
rely on more specific covariance bounds). 
\hfill $\Box$ }
\end{Remark}

\bigskip
A slightly cleaner version of the upper bound in Theorem \ref{mainres} (at the cost of slightly
worse constants) can be obtained as follows. From Schwarz inequality we see that
\begin{eqnarray*}
\lefteqn{\int \mbox{E}[|D_z F|^2 |D_z L^{-1} (F-\mbox{E}(F)) ])  d\mu(z)\leq 
\left(\int 
\mbox{E}(D_z F)^4 d\mu(z) \right)^{1/2}}\hspace*{7cm}\\&& \times \left(\int \mbox{E} |D_z L^{-1} (F-\mbox{E}(F))|^2 d\mu(z)\right)^{1/2},
\end{eqnarray*}
while, on the other hand, $\int \mbox{E} |D_z L^{-1} (F-\mbox{E}(F))|^2 d\mu(z)=\sum_{n=1}^\infty (n-1)! \|f_n\|_n^n\leq 
\sum_{n=1}^\infty n! \|f_n\|_n^2=\mbox{Var}(F)$. From this, arguing as in the last proof, we see that under the assumptions
of Theorem \ref{mainres} 
\begin{equation}\label{mainresold}
d_W\left(\frac{F_\lambda(\eta_\lambda)-\mbox{E}(F_\lambda(\eta_\lambda))}{\sqrt{\mbox{Var}(F_\lambda(\eta_\lambda))}},X \right)\leq
\frac 1{\sqrt{\lambda}} \frac{C_{2}}{C_1(\lambda)},\quad \lambda >0,
\end{equation}
with 
$C_2={C}_{2,c}^{1/2}+C_{2,b}^{1/2}$, 
\begin{equation}\label{C2c}
{C}_{2,c}=\ell(\mathcal{A})\int_{(\mathbb{R}^d)^3} {W}_c(z_1,z_2,z_3) dz_1 dz_2 dz_3,
\end{equation}
and
$${W}_c(z_1,z_2,z_3)=\int_{\mathcal{B}^4}  e^{-\ell(\cup_{i=0}^3 Q_i)}  \ell(\cap_{i=0}^3Q_i) \prod_{i=0}^3 \nu(ds_i),$$
where, as before, in the last integral, $z_0=0$.

\bigskip
We consider next the case of the Kolmogorov distance. 
Apart from the above defined $C_{2,a}$, $C_{2,b}$ and $C_{2,c}$, the relevant constants in this case are 
\begin{equation}\label{C2d}
C_{2,d}=\ell(\mathcal{A})\int_{(\mathbb{R}^d)^5} W_{d}(z_1,\ldots,z_5)dz_1\cdots dz_5,
\end{equation}
with
\begin{eqnarray*}
\lefteqn{W_{d}(z_1,\ldots,z_5)=\int_{\mathcal{B}^6} e^{-\ell( \cup_{i=0}^5 Q(z_i,s_i))} \ell(\cap_{i=0}^2 Q(z_i,s_i))
\ell(\cap_{i=3}^5 Q(z_i,s_i)) }\hspace*{3cm}\\
&&\times\mathbf{1}\big((\cup_{i=0}^2 Q(z_i,s_i) \cap (\cup_{i=3}^5 Q(z_i,s_i)))\ne \emptyset \big) \prod_{i=0}^5 \nu(ds_i),
\end{eqnarray*}
and
\begin{equation}\label{C2e}
C_{2,e}=\ell(\mathcal{A})\int_{(\mathbb{R}^d)^4} W_{e}(z_1,z_2,z_3)dz_1dz_2dz_3,
\end{equation}
with
$$W_{e}(z_1,z_2,z_3)=\int_{\mathcal{B}^4} e^{-\ell( \cup_{i=0}^3 Q(z_i,s_i))} \ell(\cap_{i=0}^2 Q(z_i,s_i))
\ell(\cap_{i=1}^3 Q(z_i,s_i))  \prod_{i=0}^3 \nu(ds_i).$$
Again $z_0=0$ in the above integrals. With this notation we have the following result.
\begin{Theorem}\label{mainresK}
Under the assumptions of Theorem \ref{mainres}, for $\lambda>0$,
$$d_K\Big(\frac{F_\lambda(\eta_\lambda)-\mbox{\em E}(F_\lambda(\eta_\lambda))}{\sqrt{\mbox{\em Var}(F_\lambda(\eta_\lambda))}},X \Big)\leq
\frac 1 {C_1(\lambda)}\frac 1{\sqrt{\lambda}} \Big({\textstyle \frac{4+\sqrt{2\pi}}{8}\frac{C_{2,a}}{C_1(\lambda)^{1/2}}+C_{2,b}^{1/2} 
+\frac 1 2 \frac{C_{2,d}^{1/2}}{\lambda^{1/2}C_1(\lambda)^{1/2}} +(C_{2,c}+9C_{2,e})^{1/2}}\Big),$$
with $C_{1}(\lambda)$ as in (\ref{C1lambda}) and constants $C_{2,a}$ to $C_{2,e}$ defined as in
(\ref{constantsC2}), (\ref{C2c}), (\ref{C2d}) and (\ref{C2e}).
\end{Theorem}

\medskip
\noindent\textbf{Proof.} We use the second inequality in Theorem \ref{WassersteinSteinPoisson}. The first two terms can be handled
as in Theorem \ref{mainres}. For the third term we use Fubini's Theorem and Lemma \ref{momentosPoisson} to get 
\begin{eqnarray*}
\lefteqn{\mbox{Var}\big({\textstyle \int_{\mathbb{R}^d}} (D_zF_\lambda(\eta_\lambda)) )^2 
|D_z L^{-1}(F_\lambda(\eta_\lambda)-\mbox{E}F_\lambda(\eta_\lambda) ) |\lambda dz\big)}\hspace*{1cm}\\
&\leq  &
\int_{(\mathcal{A}\times\mathcal{B})^6} \lambda \ell(\cap_{i=0}^2 Q_\lambda(x_i,s_i)) 
 \lambda \ell(\cap_{i=3}^5 Q_\lambda(x_i,s_i)) e^{-\lambda\ell(\cup_{i=0}^5 Q_\lambda(x_i,s_i))}\\
 &\times & \mathbf{1}\big((\cup_{i=0}^2 Q_\lambda(x_i,s_i))\cap(\cup_{i=3}^5 Q_\lambda(x_i,s_i))\ne \emptyset \big)\prod_{i=0}^5(dx_i  d\nu(s_i)).
\end{eqnarray*}
Again, the change of variables, $x_i=x_0+\lambda^{-1/d}z_i$, $i=1,\ldots,5$, yields
\begin{equation}\label{W2c}
\mbox{Var}\big({\textstyle \int_{\mathbb{R}^d}} (D_zF_\lambda(\eta_\lambda)) )^2 
|D_z L^{-1}(F_\lambda(\eta_\lambda)-\mbox{E}F_\lambda(\eta_\lambda) ) |\lambda dz\big)\leq \frac{C_{2,d}} {\lambda^5} .
\end{equation}
Finally, we turn to the fourth term. We note that the proof of Lemma \ref{momentosPoisson} shows that
\begin{eqnarray*}
\mbox{E}\left(\mathbf{1}(\eta(A)=0) \int_0^1 u^{\eta(B)}e^{-\mu(B)u}du
\mathbf{1}(\eta(C)=0) \int_0^1 v^{\eta(D)}e^{-\mu(D)v}dv \right)\leq e^{-\mu(A\cup B\cup C\cup D)}.
\end{eqnarray*}
Using this fact, Fubini's Theorem and arguing as above we see that
\begin{equation}\label{W2d1}
E\big({\textstyle \int_{\mathbb{R}^d}} 
(D_z F_\lambda(\eta_\lambda))^2 (D_zL^{-1}(F_\lambda(\eta_\lambda)-E(F_\lambda(\eta_\lambda)))^2 dz\big)
\leq  \frac{C_{2,c}}{\lambda^3}.
\end{equation}
For the double integral we use  that $D_x(UV)=(D_xU) (D_xV)+(D_xU)V+U(D_xV)$ to get that 
$D_y(C(x))=D_{x,y}^2F D^2_{x,y}L^{-1}(F-\mbox{E}F)+D_{x,y}^2F D_{x}L^{-1}(F-\mbox{E}F)+
D_{x}F D^2_{x,y}L^{-1}(F-\mbox{E}F)$ and similary for $D_y(C(x))$. The expected value of the double
integrals of the nine terms in the resulting cross-product can be handled similarly. For instance,
using the pathwise expressions for $D_{x} F_\lambda(\eta_\lambda)$ and 
$D_{x}L^{-1}(F_\lambda(\eta_\lambda)-E(F_\lambda(\eta_\lambda))$ we obtain similar pathwise expressions for
$D^2_{x,y} F_\lambda(\eta_\lambda)$ and $D^2_{x,y}L^{-1}(F_\lambda(\eta_\lambda)-E(F_\lambda(\eta_\lambda))$, namely,
$$D^2_{x,y} F_\lambda(\eta_\lambda)=\int_{\mathcal{A}\times\mathcal{B}} \mathbf{1}(\eta_\lambda (Q_\lambda(x_1,s_1)=0) \mathbf{1}(x\in Q_\lambda(x_1,s_1))
\mathbf{1}(y\in Q_\lambda(x_1,s_1))dx_1d\nu(s_1),$$
\begin{eqnarray*}
\lefteqn{D^2_{x,y}L^{-1}(F_\lambda(\eta_\lambda)-E(F_\lambda(\eta_\lambda))}\hspace*{1cm}\\
&=&-\int_{\mathcal{A}\times\mathcal{B}} \Big(\int_0^1 (1-u)u^{\eta_\lambda(Q_\lambda(x_2,s_2))}e^{-u\lambda\ell(Q_\lambda(x_2,s_2))}du\Big)\\
&&\times
\mathbf{1}(x\in Q_\lambda(x_2,s_2))
\mathbf{1}(y\in Q_\lambda(x_2,s_2))dx_2d\nu(s_2).
\end{eqnarray*}
This allows to use similar arguments as above to conclude that
\begin{eqnarray}\label{W2d2}
\lefteqn{E\big({\textstyle \int_{\mathbb{R}^d}} 
(D^2_{x,y} F_\lambda(\eta_\lambda))^2 (D^2_{x,y}L^{-1}(F_\lambda(\eta_\lambda)-E(F_\lambda(\eta_\lambda)))^2 \lambda dx \lambda dy\big)}\hspace*{2cm}\\
&\leq &  \frac{\ell(\mathcal{A})}{\lambda^3} \int_{\mathcal{A}^3\times \mathcal{B}^4} e^{-\ell(\cup_{i=0}^3 Q(z_i,s_i))} \ell(\cap_{i=0}^3 Q(z_i,s_i))^2
\prod_{i=1}^3 dz_i \prod_{i=0}^3 d\nu(s_i)\\
&\leq & \frac{C_{2,e}}{\lambda^3}.
\end{eqnarray}
An analogous approach can be used to cover the other terms in the cross product and complete the proof. \quad $\Box$

\medskip
\begin{Remark}\label{nonhomogeneouscase}{\em
While the formulation of Theorem \ref{mainres} concerns homogeneous Poisson processes
on $\mathbb{R}^d$, similar bounds hold if $\eta_\lambda$ is a Poisson point process on 
$\mathbb{R}^d$ with
a nonatomic, $\sigma$-finite intensity measure $\mu_\lambda$. Condition (\ref{fintemean}) now
reads
\begin{equation}\label{finitenhmean}
\int_{\mathcal{A}\times\mathcal{B}}e^{-\mu_\lambda(Q_{\lambda}(x,s)) }dx d\nu(s)<\infty, 
\end{equation}
while if further
\begin{equation}\label{finitenhvar}
\tilde{C}_1(\lambda):=\int_{U_\lambda} V_\lambda(x,z)dx dz<\infty,
\end{equation}
where
$$V_\lambda(x,z)=\int_{\mathcal{B}\times\mathcal{B}}
e^{-(\mu_\lambda(Q_\lambda(x,s))+\mu_\lambda(Q_\lambda(x+\lambda^{-1/d}z,t))) }
\big[ e^{\mu_\lambda(Q_\lambda(x,s)\cap Q_\lambda(x+\lambda^{-1/d}z,t))}-1 \big]d\nu(s)d\nu(t),$$
then $F_\lambda(\eta_\lambda)$ has finite second moment and an application
of Theorem \ref{WassersteinSteinPoisson} and the argument leading to (\ref{mainresold}) yields
\begin{equation}\label{otherbound}
d_W\left(\frac{F_\lambda(\eta_\lambda)-\mbox{E}(F_\lambda(\eta_\lambda))}{\sqrt{\mbox{Var}(F_\lambda(\eta_\lambda))}},X \right)\leq
\frac 1{\sqrt{\lambda}} \frac{\tilde{C}_{2,a}^{1/2}(\lambda)+
\tilde{C}_{2,b}^{1/2}(\lambda)}{\tilde{C}_1(\lambda)},\quad \lambda >0,
\end{equation}
with $\tilde{C}_{2,a}^{1/2}(\lambda)=\int_{V_\lambda}W_{\lambda,a}(x,z_1,z_2,z_3)dxdz_1dz_2dz_3$,
$$W_{\lambda,a}(x,z_1,z_2,z_3)=\int_{\mathcal{B}^4}  e^{-\mu_\lambda(\cup_{i=0}^3 Q_i)}  \mu_\lambda(\cap_{i=0}^3B_i) \prod_{i=0}^3 \nu(ds_i),$$
($Q_0=Q_\lambda(x,s_0)$, $Q_i=Q_\lambda (x+\lambda^{-1/d}z_i,s_i)$, $i=1,2,3$, in the last integral) and a similar
definition for $\tilde{C}_{2,a}^{1/2}(\lambda)$. Now, (\ref{otherbound}) becomes an interesting upper bound if 
$\frac{\tilde{C}_{2,a}^{1/2}(\lambda)+
\tilde{C}_{2,b}^{1/2}(\lambda)}{\tilde{C}_1(\lambda)}$ is a bounded function of $\lambda\in[\lambda_0,\infty)$ 
for some $\lambda_0>0$. If the intensity measures $\mu_\lambda$ are absolutely continuous
with respect to Lebesgue measure and $d\mu_\lambda(x)=\lambda f(x)dx$ for some locally integrable $f$, then,
provided $Q(s)$ is a bounded open set that contains $0$, we have by the Lebesgue differentiation theorem 
(see, e.g., Theorem 3.21 in \cite{Folland}) that for almost every $x$
$$\lim_{\lambda\to\infty}\mu_\lambda(Q_{\lambda}(x,s))=f(x)\ell(Q(s)),\quad 
\lim_{\lambda\to\infty}\mu_\lambda(Q_{\lambda}(x,s)\cup Q_\lambda(x+\lambda^{-1/d}z,t))=f(x)\ell(Q(0,s)\cup Q(z,t)).$$
Then, provided that we can exchange the order in which we take limits and integration we will
have
\begin{equation}\label{tildeC1}
\lim_{\lambda\to\infty} \tilde{C}_1(\lambda)=\int_{\mathcal{A}\times \mathbb{R}^d} \tilde{V}(x,z)dxdz :=\tilde{C}_1
\end{equation}
with
$$\tilde{V}(x,z):=\int_{\mathcal{B}\times\mathcal{B}}
e^{-f(x)(\ell(Q(0,s))+\ell(Q(z,t))) }\big[ e^{f(x)\ell (Q(0,s)\cap (Q(z,t)))}-1 \big]d\nu(s)d\nu(t).$$
Similarly, if we can exchange limits and integration for $\tilde{C}_{2,a}(\lambda)$ and 
$\tilde{C}_{2,b}(\lambda)$
we will get 
\begin{equation}\label{tildeC2a}
\lim_{\lambda\to\infty} \tilde{C}_{2,a}(\lambda)=\int_{\mathcal{A}\times (\mathbb{R}^d)} \tilde{W}_{a}(x,z_1,z_2,z_3)dxdz_1dz_2dz_3 :=\tilde{C}_{2,a}
\end{equation}
and 
\begin{equation}\label{tildeC2b}
\lim_{\lambda\to\infty} \tilde{C}_{2,b}(\lambda)=\int_{\mathcal{A}\times (\mathbb{R}^d)} \tilde{W}_b(x,z_1,z_2,z_3)dxdz_1dz_2dz_3 :=\tilde{C}_{2,b}
\end{equation}
with 
$$\tilde{W}_a(z_0,z_1,z_2,z_3)=\int_{\mathcal{B}^4}  f(z_0)e^{-f(z_0)\ell(\cup_{i=0}^3 Q_i)}  \ell(\cap_{i=0}^3Q_i) \prod_{i=0}^3 \nu(ds_i),$$
$$\tilde{W}_b(z_0,z_1,z_2,z_3)=\int_{\mathcal{B}^4} f(z_0)^2 e^{-f(z_0)\ell(\cup_{i=0}^3Q_i)} \ell(Q_0\cap Q_1) \ell(Q_2\cap Q_3) 
\mathbf{1}([Q_0\cup Q_1]\cap [Q_2\cup Q_3]\ne \emptyset) \prod_{i=0}^3\nu(ds_i),$$
where, as before, $Q_i$ is short notation for $Q(z_i,s_i)$ and $z_0=0$. Hence, provided that (\ref{tildeC1}),
(\ref{tildeC2a}) and (\ref{tildeC2b}) hold with $\tilde{C_1}\in(0,\infty)$,  $\tilde{C}_{2,a}<\infty$ and
$ \tilde{C}_{2,b}<\infty$, we have that for each $\varepsilon>0$ there exists $\lambda_0>0$ such that
\begin{equation}\label{otherboundb}
d_W\left(\frac{F_\lambda(\eta_\lambda)-\mbox{E}(F_\lambda(\eta_\lambda))}{\sqrt{\mbox{Var}(F_\lambda(\eta_\lambda))}},X \right)\leq
\frac {(1+\varepsilon)}{\sqrt{\lambda}} \frac{\tilde{C}_{2,a}^{1/2}+
\tilde{C}_{2,b}^{1/2}}{\tilde{C}_1},\quad \lambda \geq \lambda_0
\end{equation}
and we get a proper Berry-Esseen bound in this nonhomogeneous setup. Similar considerations apply to the case
of Kolmogorov distance. \quad $\Box$
}
\end{Remark}

\bigskip
We conclude this section with a discussion about the extension of Theorem \ref{mainres}
to the case of empirical measures. Again, for a cleaner presentation, we restrict ourselves
to the case of uniform empirical measures, that is, we will assume that $\mathcal{A}\subset\mathbb{R}^d$
is a bounded open set, $\{X_n\}_{n\geq 1}$ are i.i.d. uniform r.v.'s on $\mathcal{A}$ and denote
$\nu_n=\sum_{i=1}^n \delta_{X_i}$. We are interested in Berry-Esseen bounds for $F_n(\nu_n)$, where
$F_n$ is the avoidance functional of Theorem \ref{mainres}. We will show that, under suitable assumptions,
Berry-Esseen bounds for $F_n(\nu_n)$ can be obtained from Berry-Esseen bounds for $F_n(\eta_n)$ where
$\eta_n$ is a Poisson process with intensity $n$ on $\mathcal{A}$. This processes fit into the setup
of Remark \ref{nonhomogeneouscase}, having intensity measure $d\mu_n(x)=\frac{n}{\ell(\mathcal{A})} \mathbf{1}_\mathcal{A}(x) dx$.
Hence, if, for each $s$, $Q(s)$ is a bounded open set that contains 0 then for $x\in\mathcal{A}$ and
large enough $n$ we have $\mu_n(Q_n(x,s))=n\ell(Q_n(x,s))=\ell(Q(s))$,  
$\mu_n(Q_{n}(x,s)\cup Q_n(x+n^{-1/d}z,t))=\ell(Q(0,s)\cup Q(z,t))$ and similarly for the related
quantities in the integrals in $\tilde{C}_{2,a}(n)$ and $\tilde{C}_{2,b}(n)$. To keep this discussion 
simpler we make the assumption that
\begin{equation}\label{simplification}
\nu \mbox{ is a finite measure; }\quad \mbox{for some $K>0$ and all $s\in\mathcal{B}$, } Q(s)\subset B(0,K).
\end{equation}
Then, by dominated convergence we get
$$\lim_{n\to\infty} \tilde{C}_1(n)={C}_1,\quad \lim_{n\to\infty} \tilde{C}_{2,a}(n)=C_{2,a}, \mbox{ and }
\lim_{n\to\infty} \tilde{C}_{2,b}(n)=C_{2,b}$$
with $C_1$, $C_{2,a}$ and $C_{2,b}$ as in (\ref{C1}) and (\ref{constantsC2}) and, therefore, since this constants
are finite (by \ref{simplification}) we have that for some positive constant
\begin{equation}\label{poissontrunc}
d_W\left(\frac{F_n(\eta_n)-\mbox{E}(F_n(\eta_n))}{\sqrt{\mbox{Var}(F_n(\eta_n))}},X \right)\leq
\frac {C}{\sqrt{n}}, \quad n\geq 1.
\end{equation}
To derive a Berry-Esseen bound for for $F_n(\nu_n)$ from (\ref{poissontrunc}) we use the well known fact that
$\eta_n\overset d =\sum_{i=1}^{N_n}\delta_{X_i}$ if $N_n$ is a Poisson random variable with mean $n$ independent 
of the $X_i$'s (with $\eta_n$ the null measure if $N_n=0$). We introduce the constant
\begin{equation}\label{alphan}
\alpha_n:=\int_{\mathcal{A}\times \mathcal{B}} \big(1-{\textstyle \frac{\ell(Q_n(x,s)\cap \mathcal{A})}{\ell(\mathcal{A})}}\big)^n 
n {\textstyle \frac{\ell(Q_n(x,s)\cap \mathcal{A})}{\ell(\mathcal{A})} }dxd\nu(s)
\end{equation}
and note that $\alpha_n\leq \tilde{K} \nu(\mathcal{B})<\infty$, with $\tilde{K}=\ell(B(0,K))$, since we are assuming (\ref{simplification}). 
From Lemmas \ref{lematecnico3} and \ref{lematecnico4} below we see that under these assumptions there exist constants $n_0, D$ (depending only
on $K$ such that for $n\geq n_0$ and $|k-n|<n$ we have
$$\Big|\mbox{E}\big(F_n(\nu_n)-F_n(\nu_k)-\textstyle \frac{\alpha_n}n (k-n)\big) \Big|\leq D\big(\frac{k-n}{n}\big)^2$$
and 
$$\mbox{Var}(F_n(\nu_n)-F_n(\nu_k))\textstyle \leq D \frac{|k-n|}{n} \frac 1 n.$$
If we take now $\nu\in (0,\frac 1 6)$ and set $L_n=n^{1/2+\nu}$ we see that for $n\geq n_0$
\begin{eqnarray*}
\lefteqn{E\Big(\sqrt{n}(F_n(\nu_n)-F_n(\eta_n))-\textstyle \alpha_n\frac{N_n-n}{\sqrt{n}}\Big)^2=\sum_{k=0}^\infty
P(N_n=k)n E\big((F_n(\nu_n)-F_n(\nu_k))-{\textstyle \frac{\alpha_n}n} (k-n)\big)^2}\hspace*{4.5cm}\\
&=&\sum_{|k-n|\leq L_n }
P(N_n=k)n E\big((F_n(\nu_n)-F_n(\nu_k))-{\textstyle \frac{\alpha_n}n} (k-n)\big)^2\\
&+&n E\big(\big((F_n(\nu_n)-F_n(\eta_k))-{\textstyle \frac{\alpha_n}n} (N_n-n)\big)^2\mathbf{1}(|N_n-n|>L_n)\big)\\
&\leq & D\frac 1 {n^{1/2-\nu}}+D^2\frac 1 {n^{1-4\nu}}+C n (P(|N_n-n|>L_n))^{1/2}\\
&\leq & C\left(\frac 1 {n^{1/2-\nu}}+ n (P(|N_n-n|>L_n))^{1/2}\right),
\end{eqnarray*}
for some constant $C$, where the last bound comes from H\"older's inequality  and the fact that, since
$F_n(\nu_n)$ and $F_n(\eta_n)$ are bounded above by $\ell(\mathcal{A})\nu(\mathcal{B})$ and $E((N_n-n)/\sqrt{n})^4=3+\frac 1 n$,
$E\big((F_n(\nu_n)-F_n(\eta_n))-{\textstyle \frac{\alpha_n}n} (N_n-n)\big)^4$ is a bounded sequence. On the other hand, from Chernoff's
inequality we know that
$$P(|N_n-n|>L_n)\leq e^{-nh(L_n/n)}+e^{-nh(-L_n/n)}$$
with $h(u)=(1+u)\log(1+u)-u$, $u\geq -1$. A Taylor expansion shows that $h(u)\sim \frac {u^2}2$ as $u\to 0$ (the ratio tends to one)
which means that for some positive constant $c$ and large enough $n$ we have $h(L_n/n)\geq c L_n^2/n^2$, 
$h(-L_n/n)\geq c L_n^2/n^2$ and, consequently,
$$P(|N_n-n|>L_n)\leq 2 e^{-cn^{2\nu}}.$$
Hence, $n (P(|N_n-n|>L_n))^{1/2}$ vanishes at a faster rate than $n^{-(1/2-\nu)}$ and we see that for some constants (that we, again,
call $n_0$ and $C$) we have
$$E\Big(\sqrt{n}(F_n(\nu_n)-F_n(\eta_n))-\textstyle \alpha_n\frac{N_n-n}{\sqrt{n}}\Big)^2\leq \frac C {n^{1/2-\nu}},\quad n\geq n_0.$$
This last bound shows that, on one hand,
\begin{equation}\label{D1}
\lim_{n\to\infty} n \mbox{Var} (F_n(\nu_n))= C_1- \Big(\int_{\mathcal{B}} e^{-\ell(Q(s))/\ell(\mathcal{A})} \ell(Q(s))d\nu(s) \Big)^2=:D_1,
\end{equation}
with $C_1$ as in (\ref{C1}), and on the other hand
\begin{equation}\label{aproxpoisb}
d_W\big(\textstyle \sqrt{n}(F_n(\nu_n)-EF_n(\nu_n)), \sqrt{n}(F_n(\eta_n)-EF_n(\eta_n))+\textstyle 
\alpha_n\frac{N_n-n}{\sqrt{n}} \big)\leq \frac{C'}{n^{1/4-\nu/2}}.
\end{equation}
A further use of Theorem \ref{WassersteinSteinPoisson} (a straightforward modification of the argument 
in the proof of Theorem \ref{mainres}) yields that for some constants $C'',n_1$
\begin{equation}\label{aproxnormalb}
d_W\big(\sqrt{n}(F_n(\eta_n)-EF_n(\eta_n))+\textstyle 
\alpha_n\frac{N_n-n}{\sqrt{n}} ,X' \big)\leq \frac{C''}{n^{1/2}}, \quad n\geq n_1,
\end{equation}
where $X'$ is a centered normal r.v. with variance $D_1$.
Finally, combining (\ref{aproxpoisb}) and (\ref{aproxnormalb}) we conclude that
for every $\delta\in (0,\frac 1 4)$
there exists a constant $C(\delta)$ such that
$$
d_W\big(\sqrt{n}(F_n(\nu_n)-EF_n(\nu_n)) ,X' \big)\leq \frac{C(\delta)}{n^{1/4-\delta}}, \quad n\geq 1.
$$
With a simple rescaling we obtain the next result that summarizes this discussion.
\begin{Theorem}\label{mainresempir}
Assume that $\mathcal{A}\subset\mathbb{R}^d$
is a bounded open set, $\{X_n\}_{n\geq 1}$ are i.i.d. uniform r.v.'s on $\mathcal{A}$, defined on $(\Omega,\mathcal{F},P)$, 
$\nu_n=\sum_{i=1}^n \delta_{X_i}$,
$(\mathcal{B},\mathcal{G},\nu)$ a measure space and suppose further that, for each $s\in\mathcal{B}$, $Q(s)$ is a Borel set on $\mathbb{R}^d$
chosen in such a way that the map $(\omega,x,s)\mapsto \mathbf{1}(\nu_n (x+n^{-1/d}Q(s))=0)$
is $\mathcal{F}\otimes \beta^d \otimes \mathcal{G}$ measurable. If 
$$F_n(\eta)=\int_{\mathcal{A}\times\mathcal{B}} \mathbf{1}(\eta (x+n^{-1/d}Q(s))=0) d(\ell\otimes \nu)(x,s),$$
and (\ref{simplification}) holds we have (\ref{D1}) and for each $\nu\in(0,\frac 1 4)$ there exists a constant, 
$C(\delta)$ such that 
$$d_W\left(\frac{F_n(\nu_n)-\mbox{\em E}(F_n(\nu_n))}{\sqrt{\mbox{\em Var}(F_n(\nu_n))}},X \right)\leq
\frac {C(\delta)}{{n}^{1/4-\delta}},\quad n\geq 1,$$
where $X$ denotes a standard normal random variable.
\end{Theorem}

\bigskip
\begin{Remark}\label{nonoptimal}{\em
Theorem \ref{mainresempir} provides non trivial Berry-Esseen bounds for avoidance functionals of empirical measures. The assumptions
are stronger than those for the case of Poisson functionals. These assumptions can possibly be relaxed. On the other
hand, we believe that the rate in Theorem \ref{mainresempir} is not optimal and we conjecture that, under suitable assumptions, 
\begin{equation}\label{conjecture}
d_W\left(\frac{F_n(\nu_n)-\mbox{E}(F_n(\nu_n))}{\sqrt{\mbox{Var}(F_n(\nu_n))}},X \right)\leq
\frac {C}{{n}^{1/2}},\quad n\geq 1,
\end{equation}
for some constant $C$. The present approach does not yield (\ref{conjecture}) and we think that it would
be of interest to devote future research to design an approach that gives this conjectured rate for empirical
functionals. \quad $\Box$
} 
\end{Remark}

\medskip
We end the section with two technical lemmas used in the proof of Theorem \ref{mainresempir}.

\begin{Lemma}\label{lematecnico3} Under the assumptions of Theorem \ref{mainresempir} there exists a constant, $D>0$ such that for $n\geq 1 $
$$\Big|E\big(F_n(\nu_n)-F_n(\nu_k)-\textstyle \frac{\alpha_n}n (k-n)\big) \Big|\leq D\big(\frac{k-n}{n}\big)^2.$$
\end{Lemma}

\medskip
\noindent
\textbf{Proof.} Set $b_n(x,s)=n\ell(Q_n(x,s)\cap \mathcal{A})/\ell(\mathcal{A})$. We note first that $b_n(x,s)\leq \ell( Q(s))/\ell(\mathcal{A})\leq \tilde{K}
/\ell(\mathcal{A})$ ($\tilde{K}=\ell(B(0,K)$) for each $x\in\mathcal{A}$  and also that
\begin{equation}\label{bias1}
E\big(F_n(\nu_n)-F_n(\nu_k)\big)=\int_{\mathcal{A}\times\mathcal{B}}\textstyle  \Big(\big(1-\frac{b_n(x,s)}n\big)^n - \big(1-\frac{b_n(x,s)}n\big)^k\Big) dxd\nu(s).
\end{equation}
 Next, for $0\leq x\leq 1$ and 
$l>0$ we have
\begin{equation}\label{desig1}
|(1-x)^l-(1-lx)|\leq \textstyle \frac {l^2} 2 x^2. 
\end{equation}
Hence, for $k>n$ (\ref{desig1}) yields $|1-(1-b_n(x,s)/n)^{k-n}-(k-n)b_n(x,s)/n|\leq \frac 1 2 \Big( \frac {k-n}n\Big)^2 \tilde{K}^2/\ell(\mathcal{A})^2$
for all $x$  and $s$ and, since $0\leq 1-b_n(x,s)/n\leq 1$, we can use (\ref{bias1}) to obtain
$$\Big|E\big(F_n(\nu_n)-F_n(\nu_k)-\textstyle \frac{\alpha_n}n (k-n)\big) \Big|\leq \frac 1 2  \frac {\tilde{K}^2\nu(\mathcal{B})} {\ell(\mathcal{A})} 
\big(\frac{k-n}{n}\big)^2.$$
If $0\leq k\leq n$ (\ref{desig1}) implies  
\begin{equation}\label{desig2}
|(1-b_n(x,s)/n)^{n-k}-1-(k-n)b_n(x,s)/n|\leq \textstyle\frac 1 2 \frac{\tilde{K}^2}{\ell(\mathcal{A}^2)} \Big( \frac {k-n}n\Big)^2, 
\end{equation}
which, in turn, implies
\begin{equation}\label{desig3}
|(1-b_n(x,s)/n)^{n}-(1-b_n(x,s)/n)^{k}| \leq \textstyle \frac{\tilde{K}}{\ell(\mathcal{A})}\frac{|k-n|}{n}+ \frac{\tilde{K}^2}{\ell(\mathcal{A}^2)} \Big( \frac {k-n}n\Big)^2.
\end{equation}
From these last two estimates we obtain the conclusion.\quad $\Box$

\medskip
\medskip
\begin{Lemma}\label{lematecnico4}
Under the assumptions of Theorem \ref{mainresempir} there exist constants, $D>0$, $n_0$
such that for $n\geq n_0$ and $|k-n|\leq n$ we have
$$\mbox{\em Var}(F_n(\nu_n)-F_n(\nu_k))\textstyle \leq D \frac{|k-n|}{n} \frac 1 n.$$
\end{Lemma}

\medskip
\noindent
\textbf{Proof.} We assume for simplicity $\ell(\mathcal{A})=1$, $\nu(\mathcal{B})=1$, the general case following with straightforward
changes. If we set $a_n(x,s,y,t)=n \ell((Q_n(x,s)\cup Q_n(y,t))\cap \mathcal{A})$ 
and $b_n(x,s)=n\ell(Q_n(x,s)\cap \mathcal{A})$
then $a_n(x,y)\leq 2\tilde{K}$, $b_n(x)\leq \tilde{K}$ with $\tilde{K}=\ell(B(0,K)))$ and
\begin{eqnarray*}
\lefteqn{\mbox{Var} (F_n(\nu_n)-F_n(\nu_k)) =\int_{(\mathcal{A}\times\mathcal{B})^2} \Big[ \textstyle \big(1-\frac{a_n(x,s,y,t)}n\big)^n\big(
1+\big(1-\frac{a_n(x,s,y,t)}n\big)^{k-n}
-2 \big(1-\frac{b_n(x,s)}n\big)^{k-n}\big)}\hspace*{4cm}\\
&&-\textstyle\big(1-\frac{b_n(x,s)}n\big)^{n}\big(1-\frac{b_n(y,t)}n\big)^{n}\big(1+
\textstyle\big(1-\frac{b_n(x,s)}n\big)^{k-n}\big(1-\frac{b_n(y,t)}n\big)^{k-n}\\
&&-2 \textstyle\big(1-\frac{b_n(x,s)}n\big)^{k-n}\big)
\Big] dxd\nu(s)dyd\nu(t).
\end{eqnarray*}
Let us assume that $k>n$ and set 
\begin{eqnarray*}
\lefteqn{V_{n,k} =\int_{(\mathcal{A}\times\mathcal{B})^2} \Big[ \textstyle \big(1-\frac{a_n(x,s,y,t)}n\big)^n\big(1+e^{-\frac{k-n}{n}a_n(x,s,y,t)}
\big(1-\frac{(k-n)a_n^2(x,s,y,t)}{2n^2}\big)}\hspace*{2cm}\\
&&\textstyle -2 e^{-\frac{k-n}{n}b_n(x,s)}\big(1-\frac{(k-n)b_n^2(x,s)}{2n^2}\big)\big)\\
&&-\textstyle\big(1-\frac{b_n(x,s)}n\big)^{n}\big(1-\frac{b_n(y,t)}n\big)^{n}\big(1+e^{-\frac{k-n}{n}(b_n(x,s)+b_n(y,t))}
\big(1-\frac{(k-n)(b_n^2(x,s)+b_n^2(y,t))}{2n^2}\big)\\
&&-2 \textstyle e^{-\frac{k-n}{n}b_n(x,s)}\big(1-\frac{(k-n)b_n^2(x,s)}{2n^2}\big)\big)
\Big] dxd\nu(s)dyd\nu(t).
\end{eqnarray*}
Then, from Lemma \ref{lematecnico2} we see that $|\mbox{Var} (F_n(\nu_n)-F_n(\nu_k))-V_{n,k}|\leq A(2\tilde{K})(7+A(2\tilde{K}))\frac {k-n}{n}\frac 1 {n^2}$.
Similarly, if we set
\begin{eqnarray*}
\lefteqn{\tilde{V}_{n,k} =\int_{(\mathcal{A}\times\mathcal{B})^2} \Big[ \textstyle \big(1-\frac{a_n(x,s,y,t)}n\big)^n\big(1+e^{-\frac{k-n}{n}a_n(x,s,y,t)}
-2 e^{-\frac{k-n}{n}b_n(x,s)}\big)}\hspace*{1cm}\\
&&-\textstyle\big(1-\frac{b_n(x.s)}n\big)^{n}\big(1-\frac{b_n(y,t)}n\big)^{n}\big(1+e^{-\frac{k-n}{n}(b_n(x,s)+b_n(y,t))}
-2 \textstyle e^{-\frac{k-n}{n}b_n(x,s)}\big)
\Big] dxd\nu(s)dyd\nu(t),
\end{eqnarray*}
then $|V_{n,k}-\tilde{V}_{n,k}|\leq 5\tilde{K}^2\frac {k-n}{n}\frac 1 {n}$. Next, we note that $\tilde{V}_{n,k}=V_{n,k,1}+V_{n,k,2}$
with
$$V_{n,k,1}=\int_{(\mathcal{A}\times\mathcal{B})^2}\textstyle\big(1-\frac{b_n(x,t)}n\big)^{n}\big(1-\frac{b_n(y,s)}n\big)^{n} 
\big( e^{-\frac{k-n}{n}a_n(x,s,y,t)} -e^{-\frac{k-n}{n}(b_n(x,s)+b_n(y,t))}\big)dxd\nu(s)dyd\nu(t)$$
and
\begin{eqnarray*}
\lefteqn{V_{n,k,2}=\int_{(\mathcal{A}\times\mathcal{B})^2}\textstyle\Big[\big(1-\frac{a_n(x,s,y,t)}n\big)^n-
\big(1-\frac{b_n(x,s)}n\big)^{n}\big(1-\frac{b_n(y,t)}n\big)^{n} \Big]}\hspace*{3cm}\\
&&\times\big(1+e^{-\frac{k-n}{n}a_n(x,s,y,t)}
-2 e^{-\frac{k-n}{n}b_n(x,s)}\big)
dxd\nu(s)dyd\nu(t).
\end{eqnarray*}
Since
\begin{eqnarray*}
0&\leq &e^{-\frac{k-n}{n}a_n(x,s,y,t)} -e^{-\frac{k-n}{n}(b_n(x,s)+b_n(y,t))}\leq \textstyle \frac{k-n}n (b_n(x,s)+b_n(y,t)-a_n(x,s,y,t))\\
&=&\textstyle \frac{k-n}n
n \ell(Q_n(x,s))\cap Q_n(y,t)\cap \mathcal{A}),
\end{eqnarray*}
we see that for each $x,s,t$ the integrand in $V_{n,k,1}$ (as a function of $y$) vanishes outside $B(x,2K/n^{1/d})$,
a set with volume $2^d\tilde{K}/n$, and is bounded by $\frac {k-n}{n} \tilde{K}$. Hence, $V_{n,k,1}\leq 2^d\tilde{K}^2 \frac{k-n}n \frac 1 n$.
Finaly, to deal with $V_{n,k,2}$ we can use Lemma \ref{lematecnico1} to see that, provided $n> 4\tilde{K}$, 
$$\textstyle\Big|\big(1-\frac{a_n(x,s,y,t)}n\big)^n-e^{-a_n(x,s,y,t)} \Big|\leq \textstyle (A(2\tilde{K}) +2\tilde{K}^2)\frac 1 n,$$
$$\textstyle\Big|\big(1-\frac{b_n(x,s)}n\big)^n\big(1-\frac{b_n(y,t)}n\big)^n-e^{-(b_n(x,s)+b_n(y,t))} \Big|\leq \textstyle
(2 A(\tilde{K}) +\tilde{K}^2+(A(\tilde{K}) +\frac 1 2 \tilde{K}^2)^2)\frac 1 n.$$
Now, the fact that 
$$|1+e^{-\frac{k-n}{n}a_n(x,s,y,t)}
-2 e^{-\frac{k-n}{n}b_n(x,s)}|\leq 2\tilde{K} \textstyle \frac{k-n}{n},$$
and the argument used in the bound for $V_{n,k,1}$ allow us to conclude that $V_{n,k,2}\leq C \frac {k-n}n\frac 1 n$
for some constant $C$.
The case $k<n$ follows as in the proof of Lemma \ref{lematecnico3}\quad $\Box$

\section{Applications.}

In this section we show the power of Theorem \ref{mainres} through its application to
two classical models in stochastic geometry: germ-grain models and quantization.
We deal first with a particular type of germ grain model given by the union of balls of a fixed
radius around Poisson centers, which we truncate to a bounded box to keep the volume of that
union finite. More precisely, we consider the functional
$$
G_\lambda(\eta)=\ell\Big(\big(\cup_{z\in\eta} B(z,t)\big) \cap \big[-{\textstyle \frac 
{\lambda^{1/d}}2},{\textstyle \frac {\lambda^{1/d}}2}\big]^d \Big),
$$
where $\eta$ is homogeneous Poisson process of unit intensity on $\mathbb{R}^d$ and
$B(z,t)$ is the ball of radius $t>0$ around $z$. As before, we write $X$
for a standard normal random variable. With this notation we have the following.
\begin{Theorem}\label{germTheo}
With the above notation, we have
$\mbox{\em E}(G_\lambda(\eta))=\lambda(1-e^{-\omega_d t^d})$, where 
$\omega_d$ denotes the volume of the $d$-dimensional unit ball and
$$\lim_{\lambda\to\infty}{\textstyle \frac 1 \lambda}\mbox{\em Var}(G_\lambda(\eta))=e^{-2\omega_d t^d}\int_{B(0,2t)} (e^{\ell(B(0,t)\cap B(z,t))}-1) dz .$$
Furthermore, for each $\lambda_0>0$ there exists a finite constant, 
$C(d,t,\lambda_0)$, depending only on  $d$, $t$ and $\lambda_0$ such that
$$\max\Bigg(d_W\left(\frac{G_\lambda(\eta)-\mbox{\em E}(G_\lambda(\eta))}{\sqrt{\mbox{\em Var}(G_\lambda(\eta))}},X \right),
d_K\left(\frac{G_\lambda(\eta)-\mbox{\em E}(G_\lambda(\eta))}{\sqrt{\mbox{\em Var}(G_\lambda(\eta))}},X \right)\Bigg)\leq
\frac {C(d,t,\lambda_0)}{\sqrt{\lambda}},\quad \lambda \geq \lambda_0.$$
\end{Theorem}

\medskip
\noindent
\textbf{Proof.} We observe first that taking $\mathcal=(0,1)^d$, $\mathcal{B}$ a set with a single element (which we denote
$0$), $\mathcal{G}$ the (only) $\sigma$-field on $\mathcal{B}$, $\nu$ the probability measure concentrated
on $0$ and $Q(0)$ the open ball centered at $0\in\mathbb{R}^d$ with radius $t$,
we have $G_\lambda(\eta)\overset d =\lambda(1-F_\lambda(\eta_\lambda))$ with $\eta_\lambda$ a Poisson
process on $\mathbb{R}^d$ with constant intensity $\lambda$ and $F_\lambda$ as in Theorem \ref{mainres}.
From Lemma \ref{moments} we obtain $\mbox{E}(F_\lambda(\eta_\lambda))=e^{-\omega_d t^d}$
and 
$$\lambda \mbox{Var}(F_\lambda(\eta_\lambda))=e^{-2 \omega_d t^d} \int_{U_\lambda} (e^{\ell(B(0,t)\cap B(z,t))}-1)dxdz=C_1(\lambda),$$
where 
$U_\lambda=\{(x,z):\, z\in (0,1)^d, z\in \lambda^{1/d}((0,1)^d-x) \}$. We note also that $C_1(\lambda)$ grows to $C_1$ 
as $\lambda\to\infty$ with
$$C_1=e^{-2 \omega_d t^d} \int_{\mathbb{R}^d} (e^{\ell(B(0,t)\cap B(z,t))}-1)dz=
e^{-2 \omega_d t^d} \int_{B(0,2t)} (e^{\ell(B(0,t)\cap B(z,t))}-1)dz<\infty,$$
since $B(0,t)\cap B(z,t)=\emptyset$ (hence, the last integrand vanishes) if $\|z\|>2t$. Thus, we can apply (\ref{mainresold})
and it suffices to show finiteness of the constants $C_{2,c}$ and $C_{2,b}$.
Now, we have 
\begin{eqnarray*}
C_{2,c}&=&\int_{(\mathbb{R}^d)^3} e^{-\ell(B(0,t)\cup (\cup_{i=1}^3 B(z_i,t)))} \ell(B(0,t)\cap (\cap_{i=1}^3 B(z_i,t)))dz_1dz_2dz_3\\
&=&\int_{(B(0,2t))^3} e^{-\ell(B(0,t)\cup (\cup_{i=1}^3 B(z_i,t)))} \ell(B(0,t)\cap (\cap_{i=1}^3 B(z_i,t)))dz_1dz_2dz_3<\infty,
\end{eqnarray*}
since, as before, $B(0,t)\cap B(z_i,t)=\emptyset$ if $\|z_i\|>2t$. On the other hand, writing $D=\{
(z_1,z_2,z_3): (B(0,t)\cup B(z_1,t))\cap (B(0,t)\cup B(z_1,t))\ne\emptyset
\}$, we have 
\begin{eqnarray*}
\lefteqn{C_{2,b}=\int_{D} e^{-\ell(B(0,t)\cup (\cup_{i=1}^3 B(z_i,t)))} \ell(B(0,t)\cap B(z_1,t))
\ell(B(z_2,t)\cap B(z_3,t))dz_1dz_2dz_3}\hspace*{1mm}\\
&\leq&\int_{B(0,2t)\times B(0,6t)^2} e^{-\ell(B(0,t)\cup (\cup_{i=1}^3 B(z_i,t)))} \ell(B(0,t)\cap B(z_1,t))
\ell(B(z_2,t)\cap B(z_3,t))dz_1dz_2dz_3<\infty,
\end{eqnarray*}
where the last bound comes from the fact that if $\ell(B(0,t)\cap B(z_1,t))>0$ then $\|z_1\|<2t$ and, if this is the case and
$(z_1,z_2,z_3)\in D$, then $z_2$ or $z_3$ must have norm less than $4t$ and if, furthermore, $\ell(B(z_2,t)\cap B(z_3,t))>0$,
then the other point must have norm less than $6t$. Finally, we take $C(d,t,\lambda_0)=(C_{2,c}^{1/2}+C_{2,b}^{1/2})/C_1(\lambda_0)$
and the result follows for the case of the Wasserstein metric. For the result in Kolmogorov
distance we note that the constant $C_{2,d}$ is upper bounded by
$$\int_{B(0,6t)^5} e^{-\ell(B(0,t))}\ell(B(0,t)^2 dz_1\cdots dz_5<\infty,$$
with a similar bound proving finiteness of $C_{2,e}$.
.\quad $\Box$

\medskip
\begin{Remark}{\em
From the computations in the proof of Theorem \ref{germTheo} we see that $C_{2,c}\leq 8^d \omega_d^4t^{4d} e^{-\omega_d t^d}$. 
Similarly, $C_{2,b}\leq 72^d \omega_d^5 t^{5d} e^{-\omega_d t^d}$.
Now, if we take $\lambda>(2t)^d$, then for each $x\in (t/\lambda^{1/d},1-t/\lambda^{1/d})$ we have $\{x\}\times B(0,t)\subset U_\lambda$ and,
as a consequence, 
$$C_{1}(\lambda)\geq e^{-2 \omega_d t^d} (1-{\textstyle \frac{2t}{\lambda^{1/d}} })^d \int_{B(0,t)}(e^{\ell(B(0,t)\cap B(z,t))}-1)  dz\geq
 e^{-2 \omega_d t^d} (1-{\textstyle \frac{2t}{\lambda^{1/d}} })^d\omega_d t^d (e^{\omega_d(\frac t 2)^d}-1) ,$$
since, for each $z\in B(0,t)$ $B(z/2,t/2)\subset(B(0,t)\cap B(z,t))$. From these estimates we see that
for $\lambda>(2t)^d$
$$
d_W\left(\frac{G_\lambda(\eta)-\mbox{E}(G_\lambda(\eta))}{\sqrt{\mbox{Var}(G_\lambda(\eta))}},X \right)\leq
\frac {C(d,t)}{\sqrt{\lambda}},$$
with $$C(d,t)=\frac{8^{d/2}e^{-\frac{\omega_d}{2}t^d}(1+ 3\omega_d^{1/2}t^{d/2})}{e^{-2 \omega_d t^d} 
(1-{\textstyle \frac{2t}{\lambda^{1/d}} })^d\omega_d t^d (e^{\omega_d(\frac t 2)^d}-1)}.$$
Of course, the constant is not optimal, but we see how easily we can get a simple explicit upper bound in Theorem
\ref{germTheo}.
}
\end{Remark}

\bigskip
Next, we provide a Berry-Esseen bound for the volume of the union of balls centered 
around the points of a empirical measure. The proof is a simple application of Theorem \ref{mainresempir}.
\begin{Theorem}\label{germTheoempir}
With the notation of Theorem \ref{germTheoempir}, if $X_1,\ldots,X_n$ are i.i.d. uniform r.v.'s on $[-\frac 
{n^{1/d}}2,$ $ \frac {n^{1/d}}2]$ and $\nu=\sum_{i=1}^n \delta_{X_i}$, we have
$$\lim_{n\to\infty}\mbox{\em E}(G_n(\nu_n))=(1-e^{-\omega_d t^d})$$
and
$$\lim_{n\to\infty}{\textstyle \frac 1 n}\mbox{\em Var}(G_n(\nu_n))=e^{-2\omega_d t^d}\Big[\int_{B(0,2t)} 
(e^{\ell(B(0,t)\cap B(z,t))}-1) dz- \omega_d^2 t^{2d}\Big].$$
Furthermore, for each $\delta\in(0,\frac 1 4)$ there exists a finite constant, 
$C(d,t,\delta)$, depending only on  $d$, $t$ and $\delta$ such that
$$d_W\left(\frac{G_n(\nu_n)-\mbox{\em E}(G_n(\nu_n))}{\sqrt{\mbox{\em Var}(G_n(\nu_n))}},X \right)\leq
\frac {C(d,t,\delta)}{n^{1/4-\delta}},\quad n\geq 1.$$
\end{Theorem}

Theorems \ref{germTheo} and \ref{germTheoempir} provide improvements
over known results. \cite{Penrose2007} contains CLT's that include both the functionals $G_\lambda(\eta)$
and $G_n(\nu_n)$ in this paper. Yet, to our knowledge, we are giving the first Berry-Esseen boudns
for functionals of this type.

\bigskip

We turn now to the quantization error functional. 
As noted in the Introduction, quantization deals with the issue of 
approximation of a continuous measure by another measure concentrated on a finite set and this
problem has applications in a variety of fields. We focus on the case of quantization around
random Poisson centers. More precisely, we assume that
$\eta_\lambda$ is a Poisson process on $\mathbb{R}^d$ with constant intensity $\lambda$ and 
consider
\begin{equation}\label{quanterror}
H(\eta_\lambda)=\int_{[0,1]^d} \min_{z\in\eta_\lambda}\|x-z \|^p dx=\sum_{z\in \eta_\lambda}
\int_{C(z,\eta_\lambda)}\|x-z \|^p dx,
\end{equation}
where $C(z,\eta_\lambda)$ is the Voronoi cell around $z$, that is, the set of points $x\in [0,1]^d$ which are
closer to $z$ than to any other point in the support of $\eta_\lambda$. 

A key fact about the functional $H$ is that for any random point measure on $\mathbb{R}^d$, $\eta$,
\begin{eqnarray}\nonumber
H(\eta)&=& \int_{[0,1]^d}  \min_{z\in \eta}\|x-z \|^p dx\\
\nonumber &=& \int_{[0,1]^d}\left(p\int_0^\infty s^{p-1} \mathbf{1}(\min_{z\in\eta}\|x-z \|>s) ds \right) dx\\
\label{expresion1} &=& \int_{(0,1)^d\times (0,\infty)} \mathbf{1}(\eta (B(x,s))=0 )dx d\nu(s),
\end{eqnarray}
where $\nu$ is the Borel measure on $(0,\infty)$ given by $d\nu(s)=ps^{p-1}ds$.
While the last expression in (\ref{quanterror}) can be used to compute moments of $F(\eta_\lambda)$ via Campbell's
Theorem (see \cite{Yukich2008}) we see from (\ref{expresion1}) that the quantization error functional 
is, up to a scaling factor, an avoidance functional as those considered in Theorem \ref{mainres}. More precisely, 
we have that 
\begin{equation}\label{expression2}
H(\eta_\lambda)= \lambda^{-p/d} F_\lambda(\eta_\lambda)
\end{equation}
if $\mathcal{A}=(0,1)^d$, $\mathcal{B}=(0,\infty)$ equiped with the Borel $\sigma$-field, $d\nu(s)=ps^{p-1}ds$, $s\in(0,\infty)$
and $Q(s)=B(0,s)$. A simple consequence of this fact and Lemma \ref{moments} is that 
$$\mbox{E}(H(\eta_\lambda))=  \lambda^{-\frac p d} \int_0^\infty p s^{p-1} e^{-\ell(B(0,s))}ds= 
\lambda^{-\frac p d}\omega_d^{-\frac p d}\Gamma (1+ {\textstyle \frac p d}),$$
where $\omega_d$ denotes the volume of the $d$-dimensional unit ball and
$$\mbox{Var}(F(\eta_\lambda))=\lambda^{-1-\frac {2p} d}C_1(\lambda)$$
with $C_1(\lambda)=\int_{U_\lambda} V(z)dx dz$,
$U_\lambda=\{(x,z): \, x\in (0,1)^d,z\in \lambda^{1/d}((0,1)^d-x))  \}$ and
$$V(z)=\int_{(0,\infty)\times (0,\infty)} p^2 u^{p-1} v^{p-1}
e^{-(\ell(B(0, u))+ \ell(B(z, v)))} (e^{\ell(B(0, u)\cap B(z, v))} -1) du\,dv.$$
As in Lemma \ref{moments} we have that $C_1(\lambda)$ grows to $C_1=\int_{\mathbb{R}^d} G(z)dz$.
We prove next finiteness of $C_1$
\begin{Lemma}\label{momentquant}
With the above notation 
$G(z)$ is integrable on $\mathbb{R}^d$ and
$$\lim_{\lambda \to \infty}\lambda^{1+\frac {2p} d} \mbox{\em Var} (H(\eta_\lambda))= \int_{\mathbb{R}^{d}} G(z) dz.$$
\end{Lemma}

\medskip
\noindent
\textbf{Proof.}
To show integrability of $G$ (the limit follows from monotone
convergence), we observe that if $u+v<\|z\|$ then $B(0, u)\cap B(z, v)=\emptyset$ and 
$e^{\ell(B(0, u)\cap B(z, v))} -1=0$. On the other hand we always have
$e^{\ell(B(0, u)\cap B(z, v))} -1\leq e^{\min(\ell(B(0, u)), \ell(B(z, v)))}$. Therefore,
\begin{eqnarray*}
V(z)&\leq & \int_{(0,\infty)\times (0,\infty)\backslash (0,\frac{\|z\|}{2})\times(0,\frac{\|z\|}{2})} p^2 u^{p-1} v^{p-1}
e^{-\max(\ell(B(0, u)), \ell(B(z, v)))}  du\,dv\\
&=&2 \int_{\frac{\|z\|}{2}}^{\infty}p u^{p-1} e^{-\omega_d u^d} \left[ \int_0^u p v^{p-1}  dv \right]du\\
&=&2p \int_{\frac{\|z\|}{2}}^{\infty}p u^{2p-1} e^{-\omega_d u^d} du.
\end{eqnarray*}
From this bound and the coarea formula we conclude that
\begin{eqnarray*}
\int_{\mathbb{R}^d}V(z)dz&\leq & d \omega_d 2 p \int_0^\infty t^{d-1}\left[\int_{\frac t 2}^{\infty}u^{2p-1} 
e^{-\omega_d u^d} du  \right] dt\\
&=& \omega_d p 2^{d+1}  \int_0^\infty u^{2p+d-1} e^{-\omega_d u^d} du\\
&=& {\textstyle \frac p d } {2^{d+1}} { \omega_d^{-\frac {2p}d }} \int_0^\infty x^{\frac{2p} d} e^{-x} dx=
{\textstyle \frac p d } {2^{d+1}} { \omega_d^{-\frac {2p}d }}\Gamma({\textstyle \frac{2p} d +1})<\infty.
\end{eqnarray*}
\quad $\Box$

We are ready now for the main result about the Poisson quantization error.
In this case, and for the sake of brevity, we restrict ourselves to the case of Wasserstein distance
but a similar analysis would yield an equivalent result in terms of Kolmogorov's distance.

\begin{Theorem}\label{TeoremaPrincipal}
If $\eta_\lambda$ is a Poisson random measure on $\mathbb{R}^d$ with constant intensity $\lambda$ and 
$$H(\eta_\lambda)=\int_{[0,1]^d} \min_{z\in\eta_\lambda}\|x-z \|^p dx$$
then there exist positive constants $C, \lambda_0$ such that
$$d_{W}\left(\frac{H(\eta_\lambda)-\mbox{\em E}(H(\eta_\lambda))}{\sqrt{\mbox{\em Var}(H(\eta_\lambda))}},X \right)\leq \frac{C}{\sqrt{\lambda}},\quad \lambda\geq \lambda_0$$
where $X$ denotes a standard normal random variable.
\end{Theorem}

\medskip
\noindent
\textbf{Proof.} In view of 
Theorem \ref{mainres}, in its version (\ref{mainresold}) and Lemma \ref{momentquant} all we have to do is to 
prove finiteness of the constants $C_{2,c}$ and $C_{2,b}$ in (\ref{constantsC2}).
This follows from Lemmas \ref{DF4bound} and \ref{innerproductbound} below. \quad $\Box$

\bigskip
As we mentioned in the Introduction a CLT for $H(\eta_\lambda)$ and related functionals 
can be found in \cite{Yukich2008}. As for the case of germ-grain models, Theorem \ref{TeoremaPrincipal}
is, to our knowledge, the first Berry-Esseen bound for this type of functional.

We conclude with the Lemmas used in the proof of Theorem \ref{TeoremaPrincipal}.

\begin{Lemma}\label{DF4bound}
If $F_\lambda(\eta_\lambda)$ is defined as in (\ref{expression2}) and $C_{2,c}$
the related constant defined in (\ref{C2c}) then $C_{2,c}<\infty$.
\end{Lemma}

\medskip
\noindent 
\textbf{Proof.} We note that 
\begin{equation}
C_{2,c}= \int_{(\mathbb{R}^d)^3 \times (0,\infty)^4} {\textstyle \left({\displaystyle \prod_{i=1}^4} p s_i^{p-1}\right)}  e^{-\ell(\cup _{i=1}^4 B(y_i,s_i))}
 \ell\left(\cap_{i=1}^4 B(y_i,s_i)\right)
\prod_{i=1}^4 dy_i \prod_{i=1}^4 ds_i,
\end{equation}
where we fix $y_1=0$. A simple computation gives that 
$$C_1=\int_{\mathbb{R}^d} E\left[A(z) B(z)^3 \right] dz$$
with
$$A(z)=\int_0^\infty p s^{p-1} \mathbf{1}(\eta_1(B(0,s))=0) \mathbf{1}(\|z\|\leq s)ds$$ and 
$$B(z)=
\int_{\mathbb{R}^d\times (0,\infty)} p s^{p-1} \mathbf{1}(\eta_1(B(y,s))=0) \mathbf{1}(\|y-z\|\leq s)dyds.$$
We observe now that $B(z)\overset{d}= B(0)$ for all $z\in \mathbb{R}^d$ (change variable $y-z=x$ and use the fact 
that a shift of $\eta_1$ is still a Poisson process on $\mathbb{R}^d$ with constant unit intensity). Hence, by Schwarz
inequality,
\begin{equation}\label{Cbound}
 C\leq (\mbox{E}(B(0)^6))^{1/2} \int_{\mathbb{R}^d} \left(E\left[A(z)^2\right]\right)^{1/2} dz
\end{equation}
Next, we show that
$B(0)$ has finite moments of all orders. In fact, let us define
$$R:=\inf\{s>0: \ell(y\in \mathbb{R}^d:\, \eta_1(B(y,s))=0, \|y\|\leq s)=0 \}.$$
$R$ is the minimal radius $s$ such that every point in the ball $B(0,s)$ has at least a point of the 
Poisson process $\eta_1$ within $s$ distance. On the other hand
\begin{eqnarray*}
B:=B(0)&=&\int_0^\infty  p s^{p-1}  \ell(y\in \mathbb{R}^d:\, \eta_1(B(y,s))=0, \|y\|\leq s)   ds\\ 
&\leq& \int_0^R p s^{p-1} \ell(B(0,s))ds =\frac{c_d p}{p+d}R^{p+d}
\end{eqnarray*}
and, therefore, it suffices to show that $R$ has finite moments of all orders. We can prove this choosing a partition
of the surface of the $(d-1)$-dimensional unit sphere (the boundary of $B(0,1)$) into $N=N(d)$ measurable regions, 
$R_1,\ldots,R_N$,  of equal area 
and diameter less than one (a proof that this can be done, together with estimates on the minimal $N(d)$ number of regions 
needed in the partition, can be found, for instance, in Lemma 21 in \cite{FeigeSchechtman2002}; see also \cite{Leopardi2006}).
We set then $S_i=\cup_{t\in[0,1]} tR_i$ and $T_i=\inf\{t>0:\, \eta_1(t S_i)>0\}$. 
Then $S_1,\ldots,S_N$ is a partition of $B(0,1)$ into regions of equal volume and 
diameter less than one. On the other hand, $T_1,\ldots,T_N$ are i.i.d. random variables with
$\mbox{P}(T_i>t)= \mbox{P}(\eta_1(tS_i)=0)=e^{-\ell(tS_i)}=e^{-t^d/N}$. Hence,
\begin{equation}\label{maxTi}
\mbox{P}(\max_{1\leq i\leq N} T_i >t)\leq N e^{-t^d/N}
\end{equation}
But now, if we take $t\geq \max_{1\leq i\leq N} T_i$, then
for each $i$ there is a point of $\eta_1$ in $t S_i$. The ball of radius $t$ centered at that point covers 
$t S_i$ ($tS_i$ has diameter less than $t$). Thus, for every $y$ in the ball $B(0,t)$ there is a point of $\eta_1$ within
$t$ distance and, consequently, $t\geq R$ and this shows that $R\leq \max_{1\leq i\leq N} T_i$ and, combined with (\ref{maxTi}),
that
\begin{equation}\label{finitemomentsR}
E(R^q)\leq \int_0^\infty q t^{q-1} N e^{-t^d/N} dt<\infty,
\end{equation}
as claimed. Finally, to see that $\int_{\mathbb{R}^d} \left(E\left[A(z)^2\right]\right)^{1/2} dz<\infty$ and
complete the proof in view of (\ref{Cbound}), we write
$A(z)=\int_{(\|z\|,\infty)} p s^{p-1}\mathbf{1}(\eta(B(0,s))=0)ds$. Hence,
$$A(z)^2=\int_{(\|z\|,\infty)^2} p^2 (st)^{p-1} \mathbf{1}(\eta(B(0,s)\cup B(0,t))=0)dsdt$$
and
\begin{eqnarray*}
\mbox{E}(A(z)^2)&=&\int_{(\|z\|,\infty)^2} p^2 (st)^{p-1} e^{-c_d (\max(s,t))^d}dsdt\\
&=& 2 \int_{\|z\|}^\infty p t^{p-1} e^{-c_d t^d} \left(\int_{\|z\|}^t p s^{p-1} ds\right) dt\\
&\leq & 2 \int_{\|z\|}^\infty p t^{2p-1} e^{-c_d t^d} dt.
\end{eqnarray*}
Since $\int_{x}^\infty  t^{2p-1} e^{-c_d t^d} dt\approx x^{2p-d} e^{-c_d x^d}$ as $x\to\infty$ (in the sense
that the ratio tends to a positive constant; this follows from l'H\^opital's rule, for instance) we see that
$(\mbox{E}(A(z)^2))^{1/2}\leq K \|z\|^{p-d/2} e^{-\frac {c_d}{2}\|z\|^d}$ for some constant $K$ and
large enough $\|z\|$. But this shows that $\int_{\mathbb{R}^d} \left(E\left[A(z)^2\right]\right)^{1/2} dz<\infty$ 
and completes the proof.\quad $\Box$

\bigskip
Finally, we prove the last technical result of this section.
\begin{Lemma}\label{innerproductbound}
If $F_\lambda(\eta_\lambda)$ is defined as in (\ref{expression2}) and $C_{2,b}$
the related constant defined in (\ref{constantsC2}) then $C_{2,b}<\infty$.
\end{Lemma}

\medskip
\noindent 
\textbf{Proof.} We note that now, fixing $y_1=0$, writing $B_i=B(y_i,s_i)$, $i=1,\ldots,4$ and 
$$g(s_1,s_2,s_3,s_4)=\int_{(\mathbb{R}^d)^3}e^{-\ell(\cup _{i=1}^4 B_i)}\ell(B_1\cap B_2)\ell(B_3\cap B_4)
\mathbf{1}(({\textstyle B_1\cup B_2})\cap (B_3\cup B_4) \ne \emptyset)\prod_{i=2}^4 dy_i,$$
we have 
$$C_{2,b}=\int_{(0,\infty)^4} g(s_1,s_2,s_3,s_4)\prod_{i=1}^4 p s_i^{p-1}ds_i.
$$
Next, observe that $\ell(B_1\cap B_2)=0$ if $\|y_2\|>s_1+s_2$, whereas if $\|y_2\|\leq s_1+s_2$ then $\|z\|\leq s_1+2s_2$
for every $z\in B_1\cup B_2$. Similarly, $\ell(B_3\cap B_4)>0$ implies that $\|y_3-y_4\|<s_3+s_4$.  
Hence, if $\ell(B_1\cap B_2)>0$, $\ell(B_3\cap B_4)>0$ and 
$(B_1\cup B_2)\cap (B_3\cup B_4) \ne \emptyset$ then $\|y_i\|\leq 2(s_1+s_2+s_3+s_4)$, $i=2,3,4$. This shows that
$$g(s_1,s_2,s_3,s_4)\leq 2^d \ell (B(0,\sum_{i=1}^4 s_i))^3 \ell (B(0,s_1\wedge s_2))\ell (B(0,s_3\wedge s_4)) e^{-\ell(B(0,\vee_{i=1}^4 s_i ))}.$$
Thus, it suffices to prove that for positive $q_i$ ,$i=1,\ldots,4$,
$$\int_{\{0<s_1<s_2<s_3<s_4 \}}e^{-c_d s_4^d} \prod_{i=1}^4  s_i^{q_i}ds_i <\infty.$$
Indeed, from iterated integration we see that the last integral equals
$$\int_{0}^{\infty} e^{-c_d s_4^d} \frac{s_4^{q_1+q_2+q_3+q_4+3}}{(q_1+1)(q_1+q_2+2)(q_1+q_2+q_3+3)}ds_4 <\infty$$
and the result follows. \quad $\Box$

\section*{Appendix.}

The following technical results have been used in the proof of auxiliary Lemmas needed
for the proof of Theorem \ref{mainresempir}.

\begin{Lemma}\label{lematecnico1}
If $0<2x<n$ and $0<2y<n$ then
$$\Big| \big(1-{\textstyle \frac x n}\big)^n-e^{-x} \big(1-{\textstyle \frac {x^2}{2n}}\big)\Big|\leq {\textstyle \frac {A(x)}{n^2}},$$
$$\Big| \big(1-{\textstyle \frac x n}\big)^n(1-{\textstyle \frac y n}\big)^n-e^{-(x+y)} \big(1-{\textstyle \frac {x^2+y^2}{2n}}\big)\Big|
\leq {\textstyle \frac {x^2y^2}{4n^2}+\frac {A(x)}{n^2}+\frac {A(y)}{n^2}+\frac {A(x)A(y)}{n^4}},$$
where $A(x)=\frac 2 3 x^3 + \frac 1 8 x^4$.
\end{Lemma}

\medskip
\noindent
\textbf{Proof.}
We observe first that for $0<x<n$ we have
\begin{equation}\label{obs1}
0\leq e^{-x}- \big(1-{\textstyle \frac x n}\big)^n=\int_x^{-n\log(1-\frac x n)} e^{-t}dt.
\end{equation}
Now, from the series expansion $-n\log(1-\frac x n)=n\sum_{k=1}^\infty \frac 1 k (\frac x n)^k$, $0<x<n$
we see that $0\leq -n\log(1-\frac x n)-(x+\frac {x^2}{2n})\leq \frac n 3 \sum_{k=3}^{\infty} (\frac x n)^k=
\frac n 3\frac {x^3/n^3}{1-x/n}$. Hence, if $0<2x<n$ we have
\begin{equation}\label{obs2}
{\textstyle x+\frac {x^2}{2n}}\leq -n\log(1-\textstyle \frac x n)\leq {\textstyle x+\frac {x^2}{2n}+\frac 2 3 \frac {x^3}{n^2}}.
\end{equation}
Noting that $\int_x^{x+\frac{x^2}{2n}} e^{-t}dt=e^{-x}(1- e^{-\frac{x^2}{2n}})$ we obtain from (\ref{obs1}) and (\ref{obs2})
that
\begin{equation}\label{obs3}
\Big| \big(1-{\textstyle \frac x n}\big)^n-e^{-x} e^{-\frac{x^2}{2n}}\Big|\leq \textstyle \frac 2 3 \frac {x^3}{n^2},\quad \mbox{ if } 0<2x<n.
\end{equation}
We can easily check that $|e^{-x}-(1-x)|\leq \frac {x^2} 2$ for $x\geq 0$ and this entails that
$|e^{-\frac{x^2}{2n}}-(1-\frac{x^2}{2n})|\leq \frac {x^4} {8n^2}$. Combining this with (\ref{obs3}) we obtain the first
inequality in the statement and from this, trivially, we get the second.\quad $\Box$.

\medskip
\begin{Lemma}\label{lematecnico2}
If $0<2x<n$ and $0<k-n<n$ then
$$\Big| \big(1-{\textstyle \frac x n}\big)^{k-n}-e^{-\frac{k-n}{n}x} \big(1-{\textstyle \frac {k-n}{n}\frac {x^2}{2n}}\big)\Big|\leq 
{\textstyle A(x) \frac {k-n}{n}\frac 1 {n^2}},$$
where $A(x)=\frac 2 3 x^3 + \frac 1 8 x^4$.
\end{Lemma}

\medskip
\noindent
\textbf{Proof.} Apply Lema \ref{lematecnico1} to $x'=\frac {k-n}{n}x$, $n'=k-n$ and note that $A$ is increasing on $[0,\infty)$, hence,
$A(x')\leq A(x)$.\quad $\Box$

\end{document}